\documentclass[12pt]{article}
\usepackage[latin9]{inputenc}
\usepackage{amsmath}
\usepackage{amssymb}
\usepackage{esint}
\usepackage{graphicx}
\usepackage{xcolor}

\usepackage{amsthm}\usepackage{amsfonts}\usepackage{latexsym}

\begin{document}

\newtheorem{theo}{Theorem}[section]
\newtheorem{prop}[theo]{Proposition}
\newtheorem{cor}[theo]{Corollary} 
\newtheorem{lem}[theo]{Lemma}
\newtheorem{rem}[theo]{Remark}
\newtheorem{con}[theo]{Conjecture}
\newtheorem{as}[theo]{Assumption}
\newtheorem{de}[theo]{Definition}

\begin{center}
{\Large \bf 
  Antiferromagnetic Covariance Structure \\
  of  
Coulomb Chain  \\
}
\end{center}

\begin{center}    
TATYANA S. TUROVA\footnote{
Mathematical Center, University of
Lund, Box 118, Lund S-221 00, 
Sweden, \newline
tatyana@maths.lth.se, tel.:+46 46 222 8543 
\newline
On leave from IMPB Russian Academy of Science Branch of KIAM RAS}
\end{center}    

\bigskip

\begin{center}
{ \it In memory of Vadim Alexandrovich Malyshev, 1938-2022
}
\end{center}

\bigskip

\setcounter{page}{1}
\renewcommand{\theequation}{\thesection.\arabic{equation}}
\setcounter{equation}{0}

\begin{abstract}
We consider a system of particles lined up on a finite interval with Coulomb 3-dimensional interactions between close neighbours, i.e. only a few other neighbours apart. This model was introduced
by Malyshev (2015) to study the flow of
charged particles. Notably even the nearest-neighbours interactions case, the only one studied previously, was proved to exhibit multiple phase transitions
depending on the strength of the external force
when the number of particles goes to infinity.

Here we show that including 
  interactions beyond the nearest-neighbours ones, leads to  qualitatively new features. The order of covariances of distances between pairs of consecutive charges is changed when compared with the former nearest-neighbours case.
  Furthermore, 
  we discover that the covariances between inter-spacings exhibit the antiferromagnetic property, namely  they periodically change sign depending on the parity of number of spacings between
them,
  while their amplitude decays.

\end{abstract}

{\it MSC2010 subject classifications.} 82B21, 82B26, 60F05.

\section{Introduction}
\subsection{Models with Coulomb interactions}

Consider a system of $N+1$ identical particles (or charges)
located at random points of the interval $[0,1]$ 
$$\bar{Y}= (Y_0=0, Y_1, \ldots, Y_{N-1},Y_N=1).$$ 
The fixed values $Y_{0}=0$ and $Y_{N}=1$ mean that at both ends there are particles with fixed positions.
The potential energy of  the system is defined as a function of its configuration 
\begin{equation}\label{S}
\bar{y} \in {\cal S}:=\{(y_{0}, \ldots , y_{N}):  0=y_{0}<\ldots<y_{N}=1\}
\end{equation}
as follows
\begin{equation}\label{U}
U(\bar{y})=  \sum_{0 \leq j<k\leq N: \ k-j \leq K} \beta_{k-j}
V(|y_{k}-y_{j}|)  + \sum_{k=1}^{N}\int_{0}^{y_k}F_sds,
\end{equation}
where the positive function $V(|y_{k}-y_{j}|) $ represents a pair-wise  interaction between the particles at $y_{k}$ and $y_{j}$ which  depends only on the distance between them. The constant $K\geq 1$ is the range of interactions
as each particle interacts with its $2K$ neighbours (except the particles close to the ends). All constants
$\beta_m$ are positive, meaning that all the charges have the same sign, and the
function $F_s$ represents an external force at the point $s\in [0,1]$.
Next we assume that 
at  positive temperature $t>0$
the configuration vector $\bar{Y}$ 
has  Gibbs distribution with density
\begin{equation}\label{dGt}
f (\bar{y}; t, {\bar{\beta}} )=\frac{1}{Z_{ \frac{1}{t} U}(N) }  e^{-\frac{1}{t} U (\bar{y})} , \ \ \ \ \ 
\bar{y} \in {\cal S},
\end{equation}
where 
the normalizing factor (the partition function)  is
\begin{equation}\label{Z*}
Z_{\frac{1}{t}U}(N)=\int_{0
	<
	y_{1}<...<y_{N-1}<
	1}e^{-\frac{1}{t}U(\bar{y})} 
d\bar{y}.
\end{equation}
Due to the scaling properties of the function $f_{\beta, \gamma}(\bar{y}, t )$ 
\begin{equation}\label{fsca}
f_{\bar{\beta}}(\bar{y}, t )=f_{\frac{\bar{\beta}}{t}}(\bar{y}, 1), 
\end{equation}
without loss of generality we shall set  $t=1$ to study the  positive temperature case. 

Here we consider
a pair-wise Coulomb repulsive interaction 
\begin{equation}
V(x)=\frac{1}{x}, \ \ x>0.   \ \ \ \  \label{V}
\end{equation}

This model can be viewed as a special case of a general class of models in statistical physics known as Coulomb Gas or Riezs Gas (after \cite{R}), which have been studied over decades, and still remain to be in  focus of mathematical physics
(as confirmed by recent reviews  \cite{Le}  or \cite{Ch}). Models with potential (\ref{V}), also called  Coulomb 3-dim potential,   describe the so-called long-range interactions. This type of interactions  arises in many areas of physics (\cite{CDR}), however the corresponding models are  most difficult for  analysis.

The model of form (\ref{dGt}) was introduced by Malyshev in \cite{M}. It combines proper 3-dim  Coulomb potential with the  assumption that the particles are on a line. The fixed positions of particles at the end-points may be interpreted as a result of an external potential confining the particles configuration within an interval (see, e.g., \cite{Du} on related experiments in physics).

The one-dimensional Coulomb gas models form a special well-studied class
following the first works by  Lenard (\cite{L}, \cite{L1}, \cite{L2}).  In setting (\ref{dGt}) such  models have been studied mostly in  the case when $K=N$, i.e., when any pair of particles interacts. The one-dimensional case is known for  
translation-symmetry breaking phenomena  under certain conditions on the variance of the charge (\cite{AM}, \cite{AGL}). Recent works \cite{BourCLT}, \cite{Bour} provide a detailed analysis of a related circular  Riezs Gas model in dimention 1.

Malyshev \cite{M} considered model (\ref{dGt}) with only the nearest neighbours interactions
i.e., when $K=1$.
However, even in  this seemingly simple  case he
discovered a surprisingly rich picture of phase transitions in 
the ground states (formally, it is model (\ref{dGt}) at zero temperature $t=0$) for the external force $F_s$ being constant in $s$, i.e., when
\begin{equation}
F_s=F(N) \geq 0, \ \ s\in[0,1], \ \ \ \  \label{F}
\end{equation}
increases in $N$. It was shown that depending on the order of $F$ in $N$ there are a few distinct phases for the distribution of particles on the interval, changing abruptly with increase of the force from uniform to condensation of the particles at one end only. Similar complex phenomena are known for  other models of Coulomb gases (e.g., \cite{AM}, \cite{J1}, \cite{J2}).

The results of \cite{M} inspired the further study in \cite{MZ} of the nearest-neighbours interactions model (\ref{dGt}) at any 
positive temperature but without external force ($F_s\equiv 0$). The phase transitions with respect to the external force in form (\ref{F})  were described in \cite{T} at any temperature but again only for the case when $K=1$, i.e., only for the nearest-neighbours interactions. The system was proved to exhibit very small fluctuations, namely of order $N^{-3}$ for the variance of interspacings (\cite{M},  \cite{MZ}, \cite{T}). The same order  $N^{-3}$  of
variance but for a certain mean-field linear statistics of the one-dimensional jellium  model was recently reported in \cite{FM}.

Here we extend the analysis of the model  (\ref{dGt}) beyond the nearest-neighbours interactions. It turns out that by only taking into consideration the second nearest-neighbours interactions we encounter new phenomena. 
The most remarkable is the discovery of anti-ferromagnetic properties of the Coulomb ensemble. Observe that although the charges have the same sign in our system, the Coulomb interactions are repulsive. 
Notably anti-ferromagnetic property was recently confirmed experimentally for string correlations in ultracold Fermi-Hubbard chains \cite{H}.

Before we state our results we shall briefly explain in the next section why the case of  nearest-neighbours interactions, i.e., $K=1$ in (\ref{dGt}), is so drastically distinct from the ones with $K>1$.

\subsection{Representation via conditional distribution}
For any $K\geq 1$  let us arbitrarily  fix the vector $\bar {\beta}=(\beta_1, \ldots, \beta_K)$ with positive entries, and define an auxiliary  random vector $\bar{X}=(X_1, \ldots, X_N)$
with  density of distribution
\begin{equation}\label{fXg}
f_{\bar {X}}(\bar{x})=
f_{\bar {X}}(\bar{x})=\frac{1}{Z_{\bar {\beta},F}(N)}  e^{-\sum_{k=1}^{N-K+1} \sum_{j=0}^{K-1} \frac{\beta_{j}}{
		x_{k}+x_{k+1}+\ldots +x_{k+j}}-\sum_{k=1}^{N}F(x_{1}+\ldots +x_{k}) 
}
\end{equation}
in $\in [0,1]^N$, where $Z_{\bar {\beta}, F}(N)$ is the normalizing constant. 
Hence, the vector $\bar{X}=(X_1, \ldots, X_N)$  has  distribution of spacings as for the charges defined above 
but under the assumption that only at one end, say at the origin, there is a fixed charge, while the other end (at $x_{1}+\ldots+x_N$) is "free" on $\mathbb{R}^+$, instead of being fixed at $1$ as in  configuration $\bar{Y}$.
The relations between these vectors are clarified in 
the following statement
which   is straightforward to check.
\begin{prop}(\cite{MZ}, \cite{T})\label{AuP}
	The following distributional identity holds 
	for all $N \geq 2$: 
	\begin{equation}\label{cond}
	(1-Y_{N-1}, Y_{N-1}-Y_{N-2}, \ldots ,  Y_{2}-Y_{1}, Y_{1})\  \  \stackrel{d}{=}  \  \ 
	(X_1 , \ldots , X_{N}) \mid _{\sum_{i=1}^{N} X_i =1 }.
	\end{equation}
\end{prop}

\hfill$\Box$

Simplifying (\ref{fXg}) we get 
\begin{equation}\label{fX1}
f_{\bar {X}}(\bar{x})=\frac{1}{Z_{\bar {\beta}, F}(N)}  e^{-\sum_{k=1}^{N-K+1} \sum_{j=0}^{K-1} \frac{\beta_{j}}{
		x_{k}+x_{k+1}+\ldots +x_{k+j}} - \sum_{k=1}^{N} F(N-k+1)x_{k}
} .
\end{equation}
This shows that if and only if $K=1$, the components of the vector  $\bar{X}=(X_1, \ldots, X_N)$ are independent random variables. Additionally, when  $F=0$, they are also identically distributed.

Otherwise, when $K>1$ the components of the vector  $\bar{X}=(X_1, \ldots, X_N)$ are {\it dependent}, and moreover they are {\it not } identically distributed even if $F=0$.

This explains why the models with $K>1$ might have different behaviour when compared with the case $K=1$, i.e., only the nearest neighbours interactions. It is natural to 
begin with the case  $K=2$ and zero external force, i.e., when  $F_s\equiv 0$ in (\ref{U}). 
Denote for this particular case 
\begin{equation}\label{fX0}
f_{\bar {X}}(\bar{x})=\frac{1}{Z_{\beta, \gamma} }   e^{-
	\beta \sum_{k=1}^N\frac{1}{x_{k}}
	+ \gamma \sum_{k=1}^{N-1}\frac{1}{x_{k}+x_{k+1} }
} .
\end{equation}
This  form is in a direct relation with the previously studied case $K=1$ (\cite{M}, \cite{MZ}, \cite{T}) which is equivalent  to $\gamma =0$ in the last formula.

\subsection{Results}
We study particular cases of potential energy (\ref{U}) for $K=2$ and $F_s\equiv 0$. We begin with the case which allows an exact (in some sense) solution. Namely, to eliminate the boundary conditions we introduce a circular potential. In other words, we consider the charges  on a ring of length 1. 

\begin{theo}{\bf (Circular potential, zero external force.)}\label{T1}
	For any $\beta>0$, $0\leq \gamma \leq \beta$ and $N>2$ let
	\begin{equation}\label{U1}
	U_{\beta,\gamma}^{\circ}(\bar{y})=  \sum_{j=1}^{N} \frac{\beta}{y_{j}-y_{j-1}}
	+\sum_{j=2}^{N} \frac{\gamma}{y_{j}-y_{j-2}} + \frac{\gamma}{1-y_{N-1}+y_1},
	\end{equation}
	\[
	\bar{y} \in {\cal S}=\{(y_{0}, \ldots , y_{N}):  0=y_{0}<\ldots<y_{N}=1\},
	\]
	and let 
	$\bar{Y}^{\circ}$ be a random vector 
	on $ {\cal S}$ with Gibbs density  
	\[
	f_{\bar{Y}^{\circ}}(\bar{y})= \frac{1}{Z^{\circ}(N) }  e^{- U_{\beta,\gamma}^{\circ}(\bar{y})} , \ \ \ \ \ 
	\bar{y} \in {\cal S},
	\]
	where $Z^{\circ}(N)$ is the partition function.
	Denote 
	\begin{equation}\label{del}
	\delta =\delta(\gamma)=  \frac{\gamma}{4\beta +\gamma+ 2\sqrt{ 4\beta^2+2\beta \gamma}}.
	\end{equation}
	Then $\bar{Y}^{\circ}$ is evenly spaced in average so that
	\begin{equation}\label{TCM}
	{\mathbb{E}}\left(Y_{k+1}^{\circ}-Y_k^{\circ} \right) =\frac{1}{N},
	\end{equation}
	and $\bar{Y}^{\circ}$ has the following covariance structure: for all
	$j\leq k $
	\begin{equation}\label{TC1}
	{\bf Cov}\left(Y_{j+1}^{\circ}-Y_j^{\circ}, Y_{k+1}^{\circ}-Y_k^{\circ} \right) = \frac{1}{N^3} 
	\frac{1}{2\beta + \gamma (1 - \delta)/2}
	\left(
	(-\delta)^{(j-k)_{N/2}}+ O\left(N^{-\frac{1}{6}}
	\right)
	\right),
	\end{equation}
	where 
	\[(j-k)_{N/2} = \min\{j-k,N-(j-k)\}\]
	is the smallest number of charges on the circle between the $j$-th and the $k$-th ones.
	In particular,
	\begin{equation}\label{VarT}
	{\bf Var}\left(Y_{j+1}^{\circ}-Y_j^{\circ}\right) 
	= \frac{1}{N^3} 
	\frac{1}{2\beta + \gamma (1 - \delta)/2}\left(1 +O\left(N^{-\frac{1}{6}}
	\right)
	\right). 
	\end{equation}
	
\end{theo}

\bigskip

Here it is most natural to set  $\gamma=\beta$. However, keeping $\gamma$ as a free  parameter allows one to derive as a particular case $\gamma=0$ of Theorem \ref{T1} the results for model \cite{MZ} (or \cite{T})   with nearest-neighbour interactions (and zero external force). Thus setting $\gamma=0$ in (\ref{VarT}) we recover the corresponding results for the variance (the principal term) derived previously in \cite{MZ} and \cite{T}. Notice that the covariance (\ref{TC1}) was not considered previously even for $\gamma=0$. 

By definition (\ref{del}) it holds that  $0<\delta(\gamma)<1$ for all $\gamma>0$, while $\delta(0)=0$. Hence formula (\ref{TC1}) tells us that for any $\gamma >0$ 
the leading term in the covariance between spacings
$Y_{j+1}^{\circ}-Y_j^{\circ}$ and $Y_{k+1}^{\circ}-Y_k^{\circ}$
has the following order in $N$:
\begin{equation}\label{order}
\begin{array}{rll}
& N^{-3} \
\frac{1}{2\beta + \gamma (1 - \delta)/2} , & \mbox{ if } j=k \mbox{ (variance), }\\ \\
(-\delta)^{k-j}& N^{-3}\
\frac{1}{2\beta + \gamma (1 - \delta)/2}, & \mbox{ if } 0<k-j\ll \log N.
\end{array}
\end{equation}
Furthermore, we can add that in the proof of Theorem \ref{T1} we show that the correction term $O\left(N^{-\frac{1}{6}}
\right)$ in (\ref{TC1}) equals 
\[O\left(N^{-\frac{1}{6}}
\right)
-\frac{1}{N} 
\frac{1-\delta}{1+\delta} +o(1/N).\]
Our conjecture is that this correction term has order at least $N^{-1}$, which in turn would imply that the 
covariance between spacings
$Y_{j+1}^{\circ}-Y_j^{\circ}$ and $Y_{k+1}^{\circ}-Y_k^{\circ}$
does not decrease below the order $N^{-4}$ when  $k-j\gg \log N.$

Theorem \ref{T1} shows that a positive $\gamma$ affects the covariances between spacings which are separated by at most $\log N$ charges (or spacings) in a peculiar manner of changing the sign of correlations, and moreover decreasing exponentially order of the amplitude from $N^{-3}$ (at least as long as $0<k-j\ll \log N$). 

In case when $\gamma=0$, the covariances between any pair of spacings have order at most $N^{-3-\frac{1}{6}}$. Hence, adding only the next to the nearest-neighbours interactions increases the order of covariances in $N$, showing another kind of phase transition. 

The emergence of periodic in sign structure of the correlations when two interspacings are positively or negatively correlated depending on the parity of number of spacings between them, manifests antiferromagnetic property of the model. This property has not yet been reported for the Coulomb gas models. On the other hand, recent experiments \cite{H} reveal
antiferromagnetic structure of magnetic correlations via quantum gas microscopy of hole-doped ultracold Fermi-Hubbard chains. 

Our analysis does not imply any conclusion about the sign of correlations for  an arbitrary  finite range $K$ in (\ref{dGt}). However, the exponential decrease in value to a certain fixed order in $N$, seem to be universal.

The below established  approximation of the original model $\bar{Y}$ by the circular version
$\bar{Y}^{\circ}$ described in Theorem \ref{T1},
allows us to  describe the covariance structure of  $\bar{Y}$ in terms of 
the covariances of $\bar{Y}^{\circ}$.

\begin{theo}\label{T2o}
	For any $0\leq \gamma \leq \beta$
	for all 
	$N>2$ let
	\[
	U_{\beta,\gamma}(\bar{y})=  \sum_{j=1}^{N} \frac{\beta}{y_{j}-y_{j-1}}
	+\sum_{j=2}^{N} \frac{\gamma}{y_{j}-y_{j-2}} , \ \ \ \ \ 
	\bar{y} \in {\cal S},
	\]
	and let 
	$\bar{Y}$ be a random vector 
	on $ {\cal S}$ with Gibbs density  
	\[
	f_{\bar{Y}}(\bar{y})= \frac{1}{Z(N) }  e^{- U_{\beta,\gamma}(\bar{y})} , \ \ \ \ \ 
	\bar{y} \in {\cal S},
	\]
	where $Z(N)$ is the partition function.
	
	There are positive $\alpha(\gamma)$ and $C$ such that the following statements for the distribution of $\bar{Y}$ hold.
	$\bar{Y}$  is asymptotically evenly spaced in averaged so that
	\begin{equation}\label{TCMn}
	{\mathbb{E}}\left(Y_{k+1}-Y_k \right) = \frac{1}{N}\left( 1+ O\left(e^{-\alpha \min\{j, \ N-k\}}\right) \right)\left( 1+ O\left(N^{-1/3}\right) \right),
	\end{equation}
	and for all 
	$j\leq k $
	\begin{equation}\label{TC1n}
	{\bf Cov}\left(Y_{j+1}-Y_j, Y_{k+1}-Y_k \right)
	\end{equation}
	\[ = {\bf Cov}\left(Y_{j+1}^{\circ}-Y_j^{\circ}, Y_{k+1}^{\circ}-Y_k^{\circ} \right)
	\left( 1+ O\left(e^{-\alpha \min\{j, \ N-k\}}\right) \right)\left( 1+ O\left(N^{-1/3}\right) \right).
	\]
	
	\noindent
	where the covariances ${\bf Cov}\left(Y_{j+1}^{\circ}-Y_j^{\circ}, Y_{k+1}^{\circ}-Y_k^{\circ} \right)$ are described by Theorem \ref{T1}.
	
\end{theo}

For the models with constant external force $F_s=F(N)$ one may predict results as in \cite{T}, however, this will be a subject of a separate work.

We consider here the Coulomb potential following the primary goal of the model to describe distribution of charges. However, the analysis well maybe extended to the models with Riezs potentials.

\setcounter{equation}{0}
\renewcommand{\theequation}{\thesection.\arabic{equation}}
\section{Proofs}

\subsection*{Plan of the proofs}

As representation (\ref{cond}) suggests, we shall use Lagrange multiplier technique as explained in Section \ref{LM} to compute the asymptotic for the conditional 
expectation and covariances. For any $\lambda$ we introduce a random vector $\bar{X}_{\lambda}$, and its circular version $\bar{X}_{\lambda}^{\circ}$, such that $\bar{X}_{0}^{\circ}=\bar{X}^{\circ}$. To choose optimal value of ${\lambda}$ we shall study first unconditional distribution of 
$\bar{X}_{\lambda}^{\circ}$.

Section \ref{ED} extends results of \cite{T1} on exponentially fast decay of correlations between components of $\bar{X}_{\lambda}$, and on approximation of the components of $\bar{X}_{\lambda}$ by the components of its circular 
counterpart.

In Section \ref{MM} we disclose the relation between the mean value of $\bar{X}_{\lambda}$ and the minimum of the corresponding energy function.

In Section \ref{GL} Gaussian approximation is derived for the distribution of any subvector of $\bar{X}_{\lambda}^{\circ}$ consisting of at most $N^{1/7}$ consecutive components. 

Section \ref{EC} provides details on the values of expectations and covariances, making use of the results on the inverse circular Toeplitz matrix given separately in Appendix.
Based on this, the Lagrange multiplier value is settled in Section \ref{VLM}.

Section \ref{CLT} proves Central Limit Theorem  for the sum of (notably dependent) random variables $\sum_{k=1}^N{X}_{k,\lambda}^{\circ}$.

Finally, we prove the main result in Section \ref{MR}, combining derived Gaussian approximation and CLT.

\subsection{Lagrange multiplier}\label{LM}

The representation given in  Proposition \ref{AuP} via conditional distribution enables us to use 
the Lagrange multiplier method elaborated previously in a similar context in
\cite{MZ} and \cite{T} for the case $\gamma=0$, i.e., when  $X_1, \ldots, X_N$ are $i.i.d.$ random variables. 

Let us  embed  the distribution of $\bar{X}$ given by  (\ref{fX0})  into a more general class as follows.
For any 
$\lambda \in \mathbb{R}$ define a new random vector
$\bar{X}_{\lambda}=(X_{1,\lambda}, \ldots, X_{N,\lambda})$ with 
the following density of distribution:
\begin{equation}\label{fXl}
f_{\bar{X}_{\lambda}}(\bar{x})=\frac{1}{Z_{\beta, \gamma,\lambda}(N)}  e^{-H_{\lambda}(\bar{x})}, \ \ \ \ \ 
\bar{x} \in [0,1]^N,
\end{equation}
where 
\begin{equation}\label{Hlam}
H_{\lambda}(\bar{x}) = \beta \sum_{k=1}^N
	\frac{1}{x_{k}}+\gamma \sum_{k=1}^{N-1}
	\frac{1}{ x_{k}+x_{k+1} } +\lambda\sum_{k=1}^N
	x_{k}, 
\end{equation}
and
$Z_{\beta, \gamma,\lambda}(N)$ is the normalizing constant. 
This  notation is consistent with (\ref{fX0}) as
\[f_{\bar{X}_{0}}(\bar{x})=f_{\bar{X}}(\bar{x}),\]
i.e., $$
\bar{X}\stackrel{d}{=}  \bar{X}_{0}.$$

The remarkable property of new random vector $\bar{X}_{\lambda}$  is that 
for any $\lambda \in \mathbb{R}$
the following equality in distribution holds:
\begin{equation}\label{condl}
\bar{X} \mid _{\sum_{i=1}^{N} X_i =1 } \  \  \stackrel{d}{=}  \  \  \bar{X}_{\lambda} \mid _{\sum_{i=1}^{N} X_{i,\lambda} =1 }.
\end{equation}

Diaconis and Freedman \cite{DF} used this property to prove conditional Central Limit Theorem 
for the identically distributed random variables. 
Their result was used in \cite{MZ} for the analysis of  a Coulomb gas model. Then in \cite{T} the same technique was proved to work
even without the assumption of the identity of  distributions.
Here we further develop the argument for the case of  dependent random variables to study limit distribution. Our derivation of the Gaussian approximation is, however, different from the mentioned above approaches; it is inspired by treatise \cite{F}.

We begin with choosing parameter $\lambda$
so that the condition $
\sum_{k=1}^N X_{k,\lambda}
=1$ holds unbiased, i.e., 
\begin{equation}\label{ES}
\mathbb{E} S_{N, \lambda} : =\sum_{k=1}^N\mathbb{E} X_{k,\lambda} = 1.
\end{equation}
It turns out as we will see below, that this choice of parameter $\lambda$
facilitates analysis of the corresponding conditional distribution as it was in the previous works \cite{MZ} or \cite{T}.

The existence of $\lambda$ satisfying (\ref{ES}) is  proved below. Here we only note that by
the  properties of Gibbs distribution it holds that (straightforward to check) 
\begin{equation}\label{var}
\frac{\partial}{\partial \lambda} \mathbb{E} S_{N, \lambda}= - Var \left(S_{N, \lambda}\right).
\end{equation}
This proves that
$\mathbb{E} S_{N, \lambda} $ as a function of $\lambda$  is  strictly decreasing, hence the uniqueness of the solution $\lambda$ to (\ref{ES}) follows (when exists).
After we prove below  existence of $\lambda$ satisfying (\ref{ES}) 
we shall establish the Gaussian approximation for the distribution of $\bar{X}_{\lambda}$ 
for this value $\lambda$.
This will allow us to tackle the original distribution of $\bar{Y}$.

\subsection{Exponential decay of correlations\\
  and 
  approximation by a circular distribution}\label{ED}

As we observed already the entries of $\bar{X}_{\lambda}$ are not independent and moreover are not identically distributed, which complicates analysis.
Recent results \cite{T1} help us to use 
an approximation of the distribution of  $\bar{X}_{\lambda}$  by 
an easier ``circulant'' form of  distribution (\ref{fXl}). Namely, introduce a circulant version of function  in (\ref{Hlam}):
\begin{equation}\label{H0c}
  H_{\lambda}^{\circ, N}(\bar{x}):=H_{\lambda}(\bar{x})+\gamma \frac{1}{ x_{1}+x_{N} }
\end{equation}
  \[=\beta \sum_{k=1}^N
\frac{1}{x_{k}}+\gamma \sum_{k=1}^{N-1}
\frac{1}{ x_{k}+x_{k+1} } +\gamma 
\frac{1}{ x_{1}+x_{N} } + \lambda\sum_{k=1}^Nx_{k}, \]
\[\bar{x}=(x_1, \ldots, x_N) \in [0,1]^N,\]
and
define
$\bar{X}^{\circ}_{\lambda}$ to be a random vector whose distribution density is
\begin{equation}\label{Mt7}
f_{\bar{X}^{\circ}_{\lambda}}(\bar{x})=\frac{1}{Z^{\circ}_{\beta, \gamma,\lambda}(N)}  e^{-H_{\lambda}^{\circ, N}(\bar{x})}, \ \ \ \ \ 
\bar{x} \in [0,1]^N,
\end{equation}
where $Z^{\circ}_{\beta, \gamma, \lambda}(N)$ is the normalizing constant.

Notice that by the same argument as in Proposition \ref{AuP}  vector $\bar{X}^{\circ}_{0}=\bar{X}^{\circ}$ is related to $\bar{Y}^{\circ}$ in the same manner as explained in (\ref{cond}).

\begin{rem}\label{Ra1}
	For the ease of notation we write here $\bar{X}^{\circ}_{\lambda}=\left(X^{\circ}_{1,\lambda},\ldots X^{\circ}_{j,\lambda} , \ldots ,X^{\circ}_{N,\lambda} \right)$ and
	$H_{\lambda}^{\circ, N}=H_{\lambda}^{\circ}$. However, it should be taken into account that the distribution of $X^{\circ}_{j,\lambda}$ for any $j$ depends on $N$. When this dependence is in focus we shall write
	\[\bar{X}^{\circ , N}_{\lambda}=\left(X^{\circ ,N}_{1,\lambda}, \ldots ,X^{\circ, N}_{N,\lambda} \right).\]
	\end{rem}

        The exponential decay of dependence between the components of $\bar{X} $ (i.e., for the case $\lambda=0$ when $\bar{X}^{\circ}= \bar{X}^{\circ}_{0}$)
was established 
        in \cite{T1} in the following form. 
\begin{lem}\label{ThT10}({\protect{\cite{T1}}})  For any  $\gamma \geq  0$
  there is a positive constant $\alpha=\alpha(\gamma)$ such that for any $\beta>0$
and some 
  positive constant $C$ 
\begin{equation}\label{ed1}
\left|f_{X_{j} , X_{k}}(x,y)-f_{X_{j}}(x)f_{X_{k} }(y)\right|
\leq Ce^{-\alpha (k-j)} f_{X_{j}, X_{k} }(x,y),
\end{equation}
 for all $1\leq j<k \leq  N$ and all $(x,y) \in (0,1]^2$.
        \end{lem}
    Setting $j=1$ and $1<k\leq N/2$ in (\ref{ed1}) we get a decay of dependence  on the boundary conditions (i.e., $X_1$) as $k$ is increasing. This argument was used
         in \cite{T1} to derive from (\ref{ed1}) the following statements on
       the approximation by the circular distribution (still  for the case $\lambda=0$). 

\begin{lem}\label{Cor0}({\protect{\cite{T1}}})  For any  $\gamma \geq  0$
	there is a positive constant $\alpha=\alpha(\gamma)$ such that for any $\beta>0$ and some positive constant $C$ 
\begin{equation}\label{CM2}
\left|f_{X_{j} ,  X_{k}}(x,y)-f_{X^{\circ}_{j}, \ X^{\circ}_{k} }(x,y)\right|
\leq Ce^{-\alpha \min\{j,N-k\}} f_{X^{\circ}_{j},  X^{\circ}_{k} }(x,y),
\end{equation}
as well as  for all $N'>N$
\begin{equation}\label{CMa}
\left|f_{X^{\circ,N}_{j},\  X^{\circ,N}_{k}}(x,y)-f_{X^{\circ,N'}_{j}, \ X^{\circ,N'}_{k} }(x,y)\right|
\leq Ce^{-\alpha (j+N-k)} f_{X^{\circ,N}_{j},\  X^{\circ, N}_{k} }(x,y)
\end{equation}
for all $1\leq j<k \leq  N$ and all $(x,y) \in (0,1]^2$.
\end{lem}
Note that the second statement (\ref{CMa}) tells us about the exponential convergence
of $ f_{X^{\circ,N}_{j},\  X^{\circ, N}_{k} }$
in $N$.

Furthermore, precisely by the same arguments the following generalization of Lemma \ref{Cor0} holds as well, so we state it without a proof.
\begin{lem}\label{Cor0g}  For any  $\gamma \geq  0$
  there is a positive constant $\alpha=\alpha(\gamma)$ such that for any $\beta>0$ and some positive constant $C$ one has
	\begin{equation}\label{CM2g}
	\left|f_{X_{j} ,X_{j+1}, \ldots , \ X_{j+k}}(x_j,x_{j+1} \ldots, x_{j+k})-f_{X^{\circ}_{j},X^{\circ}_{j+1}, \ldots ,  X^{\circ}_{j+k} }(x_j,x_{j+1} \ldots, x_{j+k})\right|
	\end{equation}
	\[\leq Ce^{-\alpha \min\{j,N-j-k\}} f_{X^{\circ}_{j},X^{\circ}_{j+1}, \ldots ,\  X^{\circ}_{j+k} }(x_j,\ldots, x_{j+k}),\]
	as well as  for all $N'>N$
	\begin{equation*}\label{CMag}
	\left|f_{X^{\circ,N}_{j}, \ldots ,  X^{\circ,N}_{j+k}}(x_j, \ldots, x_{j+k})-f_{X^{\circ,N'}_{j}, \ldots ,  X^{\circ,N'}_{j+k} }(x_j, \ldots, x_{j+k})\right|
	\leq Ce^{-\alpha (N-k)} f_{X^{\circ,N}_{j}, \ldots  , X^{\circ, N}_{k} }(x_j, \ldots, x_k)
      \end{equation*}
for all $1\leq j<j+k \leq  N$ and all $x_j,x_{j+1} \ldots, x_{j+k} \in (0,1]$.
\end{lem}
\hfill$\Box$

Notice that in the above lemmas $\alpha(0)=+\infty$ since in the case $\gamma=0$ random variables
$X^{N}_{k}, k=1, \ldots, N,$ are $i.i.d.$, and also 
$\bar{X }^N  \  \stackrel{d}{=}  \  \  \bar{X}^{\circ, N} $.

Consider now $\bar{X }_{\lambda} $ for arbitrary $\lambda \geq 0$.  Observe that for any $\lambda$
\begin{equation}\label{fla1}
  f_{X_{j,\lambda} , X_{k,\lambda}}(x,y)=\frac{e^{-\lambda x - \lambda y}
    f_{X_{j} , X_{k}}(x,y)}{ {\mathbb{E}e^{-\lambda X_{j} - \lambda X_{k}}}},
\end{equation}
where $X_k, k\geq 1,$ are uniformly bounded random variables; and similar for the circular version
\begin{equation}\label{flac}
f_{X^{\circ}_{j,\lambda} , X^{\circ}_{k,\lambda}}(x,y)=\frac{e^{-\lambda x - \lambda y}
	f_{X^{\circ}_{j} , X^{\circ}_{k}}(x,y)}{ {\mathbb{E}e^{-\lambda X^{\circ}_{j} - \lambda X^{\circ}_{k}}}}.
\end{equation}
This together with
Lemma \ref{Cor0}, which says
\begin{equation}\label{CM2n}
f_{X_{j} , \ X_{k}}(x,y) =\left( 1+ O\left(e^{-\alpha \min\{j,N-k\}} \right) \right) f_{X^{\circ}_{j},\  X^{\circ}_{k} }(x,y),
\end{equation}
gives us
\begin{equation}\label{CM2n1}
f_{X_{j,\lambda} , X_{k,\lambda}}(x,y)=\left( 1+ O\left(e^{-\alpha \min\{j,N-k\}} \right) \right)
f_{X^{\circ}_{j,\lambda}, X^{\circ}_{k,\lambda} }(x,y)
\frac{{\mathbb{E}e^{-\lambda X^{\circ}_{j} - \lambda X^{\circ}_{k}}}}{{\mathbb{E}e^{-\lambda X_{j} - \lambda X_{k}}}}.
\end{equation}
Applying once again the result (\ref{CM2}) of Lemma \ref{Cor0} to the last fraction of the expectations we derive in turn an extension of the result of Lemma \ref{Cor0} for arbitrary $\lambda$:
\begin{equation}\label{CM2a}
f_{X_{j,\lambda} , X_{k,\lambda}}(x,y)=\left( 1+ O\left(e^{-\alpha \min\{j,N-k\}} \right) \right)
 f_{X^{\circ}_{j,\lambda}, X^{\circ}_{k,\lambda} }(x,y),
\end{equation}
where the correction term is uniform in $\lambda$ and $x,y.$

In the same manner we extend the result (\ref{CMa}) as well: 
for all $0<j< k <N'\leq N$
\begin{equation}\label{CMan}
f_{X^{\circ,N}_{j,\lambda},\  X^{\circ,N}_{k,\lambda}}(x,y)
=\left( 1+ O\left(e^{-\alpha (j+N'-k)} \right) \right)
 f_{X^{\circ,N'}_{j,\lambda},\  X^{\circ, N'}_{k,\lambda} }(x,y),
\end{equation}
where the correction term is uniform in $\lambda$ and $x,y.$
Combination of (\ref{CM2a}) and (\ref{CMan}) yields for all $0<j< k <N'\leq N$
\begin{equation}\label{S1}
f_{X^{N}_{j,\lambda} , X^{N}_{k,\lambda}}(x,y)=\left( 1+ O\left(e^{-\alpha (j+N'-k)} \right) \right)
f_{X^{\circ,N'}_{j,\lambda},\  X^{\circ, N'}_{k,\lambda} }(x,y),
\end{equation}
where the correction term is uniform in $\lambda$ and $x,y.$

Making use of Lemma \ref{Cor0g} and following precisely the above lines  we also derive
for all $0<j< j+k <N'\leq N$ that
\begin{equation}\label{S11}
  f_{X^{N}_{j,\lambda} , \ldots, X^{N}_{j+k,\lambda}}(x_j,\ldots, x_{j+k})=\left( 1+
    O\left(e^{-\alpha (N'-k)} \right) \right)
  f_{X^{\circ,N'}_{j,\lambda},\ldots,   X^{\circ, N'}_{j+k,\lambda}
  }(x_j,\ldots,x_{j+k}),
\end{equation}
where the correction term is uniform in $\lambda$ and $x_j, \ldots, x_{j+k} \in (0, 1].$

As an immediate corollary of the result (\ref{CM2a})
we get the following bounds 
\begin{equation}\label{Ea1}
\left| 
\mathbb{E} X^{N}_{r,\lambda} -\mathbb{E} X^{\circ, N}_{r,\lambda} 
\right|\leq Ce^{-\alpha \min\{r,N-r\}}
\mathbb{E} X^{\circ, N}_{r,\lambda},
\end{equation}
and for $j<k$
\begin{equation}\label{Ea2}
  \left|   {\bf Cov }\left({X^{N}_{j, \lambda} , X^{N}_{k, \lambda}}\right)-
    {\bf Cov }\left( {X^{\circ , N}_{j, \lambda}, X^{\circ , N}_{k, \lambda} }\right)
\right|\leq C
e^{-\alpha \min\{j,N-k\}} \left| {\bf Cov }
\left(
  {X^{\circ , N}_{j, \lambda}, X^{\circ , N}_{k, \lambda} }
    \right)\right|,
\end{equation}
where positive constant $C$ does not depend on $N$ or $\lambda$.

Next we establish also exponential decay of the involved correlations.

Making use of Lemma 
\ref{ThT10} we derive from (\ref{fla1})
the following relation for a two-points density:  for $j<k$
\begin{equation}\label{fla2}
  f_{X_{j,\lambda} , X_{k,\lambda}}(x_j, x_k)=
  \frac{e^{-\lambda x_j - \lambda x_k}
}{ {\mathbb{E}e^{-\lambda X_{j} - \lambda X_{k}}}}\left( 1+ O\left(e^{-\alpha (k-j)} \right) \right)f_{X_{j}}(x_j)  f_{X_{k}}(x_k)
\end{equation}
\[= \frac{ \mathbb{E} e^{-\lambda X_{j}} \mathbb{E} e^{- \lambda X_{k}}
}{ {\mathbb{E}e^{-\lambda X_{j} - \lambda X_{k}}}}
\left( 1+ O\left(e^{-\alpha (k-j)} \right) \right)f_{X_{j,\lambda}}(x_j)  f_{X_{k,\lambda}}(x_k).\]
Again  with a help of Lemma 
\ref{ThT10} we obtain
 \[
\frac{ \mathbb{E} e^{-\lambda X_{j}} \mathbb{E} e^{- \lambda X_{k}}
}{ {\mathbb{E}e^{-\lambda X_{j} - \lambda X_{k}}}} = 1+ O\left(e^{-\alpha (k-j)} \right), \]
and substitution of the last formula into (\ref{fla2}) gives us 
(for all large $k-j>0$)
\begin{equation}\label{Nov2}
f_{X_{j,\lambda} , X_{k,\lambda}}(x_j, x_k)=
\left( 1+ O\left(e^{-\alpha (k-j)} \right) \right)
f_{X_{j,\lambda}}(x_j)  f_{X_{k,\lambda}}(x_k).
\end{equation}

The last formula yields, in particular, an exponential decay of correlations
\begin{equation}\label{Exjk1}
{\bf Cov }\left( {X^{ N}_{j, \lambda}, X^{ N}_{k, \lambda} }\right)=
  O\left(e^{-\alpha (k-j)} \right)\mathbb{E} X_{j, \lambda} \mathbb{E} X_{k, \lambda},
\end{equation}
which together with (\ref{Ea2}) yields a corresponding bound for the circular case
\begin{equation}\label{Ea21}
 {\bf Cov }\left( {X^{\circ , N}_{j, \lambda}, X^{\circ , N}_{k, \lambda} }\right)=
 O\left( e^{-\alpha \min\{ k-j,  N-(k-j)} \} \right)
 \mathbb{E} X_{j, \lambda}^{\circ , N} \mathbb{E} X_{k, \lambda}^{\circ , N}.
\end{equation}

With a help of Lemma  \ref{Cor0g} we can extend result (\ref{Nov2})  for any finite-dimensional density
as follows.
Let $I=\{i_1, \ldots, i_{|J|}\}$ and $J=\{j_1, \ldots, j_{|J|}\}$
be two disjoint subsets of indices in 
$\{1, \ldots, N\}$ in increasing order with $i_{|I|}<j_{|J|}$.
Define a distance between these sets on the circle of $N$ vertices:
\begin{equation}\label{Nov4}
d(I,J)= \min\{i_1+N-j_{|J|}, j_{|J|}-i_{|I|}\}.
\end{equation}
Then a multidimensional analogue of (\ref{Nov2})  is given by
\begin{equation}\label{Nov3}
f_{X_{I,\lambda} , X_{J,\lambda}}(x_I, x_J)=
\left( 1+ O\left(e^{-\alpha d(I,J)} \right) \right)
f_{X_{I,\lambda}}(x_I)  f_{X_{J,\lambda}}(x_J),
\end{equation}
where  $x_{I}=(x_{i}, i\in I)$
for any set $I$.

\subsection{Minima of $H_{\lambda}$ and $H_{\lambda}^{\circ}$}\label{MM}
Due to the symmetry of function $H^{\circ}_{\lambda}$
the entries of $\bar{X}^{\circ , N}_{\lambda}$ are identically distributed (dependent) random variables. Hence,
\begin{equation}\label{Ma21}
\mathbb{E} X^{\circ, N}_{k,\lambda} =m^{\circ}_{\lambda}(N)
\end{equation}
for some positive $m^{\circ}_{\lambda}(N)$,
and a combination of (\ref{Ea1})  and
(\ref{Ma21}) gives us
\begin{equation}\label{De22n5}
\mathbb{E} X_{r,\lambda} ^{N}
=\mathbb{E} X_{N-r,\lambda} ^{N}
=\left(1+ O\left(e^{-\alpha r}\right)\right)m^{\circ}_{\lambda}(N)
\end{equation}
for $0<r\leq \frac{N}{2}$.
Therefore condition (\ref{ES}) on $\lambda$ requires roughly speaking,
\begin{equation}\label{n5}
\mathbb{E} X^{\circ, N}_{r,\lambda} =m^{\circ}_{\lambda}(N) \sim \frac{1}{N}.
\end{equation}
This instructs us to search for the solution $\lambda $ to (\ref{ES}) among $\lambda = \lambda(N) $ such that  $\lambda(N) \rightarrow \infty$. (One may as well ignore this remark as the solution  to (\ref{ES}) will be found indeed in  this range, and as we argued above, it is unique.)
It follows by the Gibbs distribution (\ref{fXl}) that for large $\lambda$ the probability measure is mostly concentrated around the minima of $H_{\lambda}$. This will be stated more precise below, but 
first we shall study the minima of $H_{\lambda}$ and $H_{\lambda}^{\circ}$.

Let us fix $\lambda$ arbitrarily. Functions $H_{\lambda}(\bar{x})$ and $H_{\lambda}^{\circ}(\bar{x})$ are both convex on
$ \mathbb{R}^N_+$  and therefore each function has its unique minimum. Denote, correspondingly, these minima 
$\bar{a}=(a_1, \ldots, a_N)$ 
 and $\bar{a}^{\circ}=(a_1^{\circ}, \ldots, a_N^{\circ})$, so that 
\begin{equation}\label{H1}
H_{\lambda}(\bar{a})=\min_{\bar{x}}H_{\lambda}(\bar{x}),
 \end{equation}
and 
\begin{equation}\label{H1cir}
H_{\lambda}^{\circ}(\bar{a}^{\circ})=\min_{\bar{x}}H_{\lambda}^{\circ}(\bar{x}).
 \end{equation}
Observe that  components 
$a_k=a_k(N)$ and $a_k^{\circ}=a_k^{\circ}(N)$
are  (in general) also functions of $N$.

By definition the vector $\bar{a}$ satisfies the system of equations
\begin{equation}\label{H3}
\frac{\partial}{\partial x_k} 
H_{\lambda}(\bar{x})\mid _{\bar{x}=\bar{a}}=0, \ \ k=1, \dots, N,
 \end{equation}
which is
\begin{equation}\label{H2}
\left\{
\begin{array}{rl}
\frac{\beta}{a_{1}^2}+ 
	\frac{\gamma}{ (a_{1}+a_{2})^2 } =\lambda, & \\ \\
\frac{\gamma}{ (a_{k-1}+a_{k})^2 } + \frac{\beta}{a_{k}^2}+ \frac{\gamma}{ (a_{k}+a_{k+1})^2 }=\lambda, & k=2, \dots, N-1,\\  \\
\frac{\beta}{a_{N}^2}+ 
	\frac{\gamma}{ (a_{N-1}+a_{N})^2 } =\lambda,& 
\end{array}
\right.
 \end{equation}
 where the first and the last equations represent the boundary conditions. 

Correspondingly, it holds that $\bar{a}^{\circ}$ satisfies 
\begin{equation}\label{H2cir}
\frac{\gamma}{ (a_{k-1}^{\circ}+a_{k}^{\circ})^2 } + \frac{\beta}{(a_{k}^{\circ})^2}+ \frac{\gamma}{ (a_{k}^{\circ}+a_{k+1}^{\circ})^2 }=\lambda,  \ \  k=1, \dots, N,
\end{equation}
where we denote  $a_0^{\circ}= a_{N}^{\circ}$ and $a_{N+1}^{\circ}= a_{1}^{\circ}$. One solution to the last system is obvious:
\[
a_k^{\circ}=\sqrt{\frac{2\beta + \gamma}{2\lambda}}=:a,\ \  k=1, \dots, N,
\]
moreover due to the uniqueness this is the only one which satisfies (\ref{H2cir}).

\begin{rem}\label{Ra2}
For a fixed $\lambda$ the  components
 $a_k^{\circ}=a_k^{\circ}(N, \lambda)=a(\lambda), k=1, \ldots, N,$ 
are constant with respect to $k$ and $N$, unlike the values $a_k(N, \lambda), k=1, \ldots, N,$ which depend also both on $k$ and $N$.
\end{rem}
 
It is natural to predict that 
\[a_k(N) \rightarrow a\]
as $k, N-k, N\rightarrow \infty$. We shall prove this later, but first we establish that if this convergence takes place it must be exponentially fast.

\begin{lem}\label{PA1}
	Suppose that for any $\varepsilon >0$ there is large $N$ and $K<N/2$ such that for all $K\leq k\leq N/2$
	\begin{equation}\label{As1}
	|a_k(N)- a|<\varepsilon,
	\end{equation}
	where 
	\begin{equation}\label{defa}
a=\sqrt{\frac{2\beta + \gamma}{2\lambda}}.
\end{equation}
	Then 
when  both $k$ and $N\rightarrow \infty$  and $k\leq N/2$, we have exponentially fast convergence:
\begin{equation}\label{H5}
|a_k- a|\leq C{\eta }^k
\end{equation}
where C is some positive constant and 
\[0< \eta = 1+4\frac{\beta}{\gamma} -\sqrt{16\left(\frac{\beta}{\gamma} \right)^2 +8\frac{\beta}{\gamma}}<1.\]
\end{lem}

\noindent
{\bf Proof.}
Denote
\begin{equation}\label{phi4}
\varphi (x)=\frac{1}{x^2}, \ \ \varphi^{-1} (x)=\frac{1}{\sqrt{x}}. 
\end{equation}

Let $\lambda$ be fixed arbitrarily. Then given $X_1>0$ let us define recurrently (when possible) positive numbers $X_2, X_3, \ldots, $ by the following system
 associated with (\ref{H2}):
\begin{equation}\label{F1}
\left\{
\begin{array}{rl}
\beta\varphi (X_1) + \gamma \varphi (X_1 +X_2)= \lambda, & \\ \\
\gamma\varphi (X_{k-1} + X_k) + \beta\varphi (X_k) +\gamma\varphi (X_k +X_{k+1})=\lambda, & k\geq 2.
\end{array}
\right.
 \end{equation}
 Observe, that for each $N$ and $\lambda$ since the vector $\bar{a}= (a_1, \ldots  , a_N)$ satisfies system (\ref{H2}), 
 for $X_1=a_1$ we have well defined $X_i=a_i$ at least for all $i< N/2$.
 
Define also
\begin{equation}\label{H21}
Y_k= \gamma\varphi (X_k +X_{k+1}), \ \ \ k\geq 1,
 \end{equation}
which for $Y_k > 0$ is equivalent to 
\[X_{k+1}= \varphi^{-1} \left(\frac{Y_k}{\gamma}\right)-X_k.\]
This notation allows us to transform (\ref{F1}) into the following dynamical system on $\mathbb{R}_+^2$
\begin{equation}\label{F2}
\left\{
\begin{array}{ll}
 X_{k}& = -X_{k-1} + \varphi^{-1} \left(\frac{Y_{k-1}}{\gamma}\right), \\ \\
Y_k & = - \beta\varphi (X_k) - Y_{k-1} +\lambda, \ \ \ \ k\geq 2, 
\end{array}
\right.
 \end{equation}
assuming that the initial condition $(X_1, Y_1)$ satisfies
$$\beta \phi(X_1)+Y_1=\lambda$$ as well as yields $(X_k, Y_k)\in \mathbb{R}_+^2$.
Then we have
\begin{equation}\label{F3}
\left\{
\begin{array}{ll}
 X_{k+1}-X_{k}& =-X_{k} + \varphi^{-1} \left(\frac{Y_k}{\gamma}\right)-\left(  -X_{k-1} + \varphi^{-1} \left(\frac{Y_{k-1}}{\gamma}\right)\right), \\ \\
Y_{k+1}-Y_k & = - \beta \varphi (X_{k+1}) - Y_{k} -\left( - \beta\varphi (X_k) - Y_{k-1} \right).
\end{array}
\right.
 \end{equation}
Writing 
\begin{equation}\label{F4}
\begin{array}{ll}
 \Delta X_{k}& = X_{k} -  X_{k-1} , \\
\Delta Y_k & = Y_k - Y_{k-1}, 
\end{array}
 \end{equation}
we get from (\ref{F3})
\begin{equation}\label{F5}
\left\{
\begin{array}{ll}
 \Delta X_{k+1}& =-\Delta X_{k} + \varphi^{-1} \left(\frac{Y_k}{\gamma}\right)-\varphi^{-1} \left(\frac{Y_{k-1}}{\gamma}\right), \\ \\
\Delta Y_{k+1} & = - \Delta Y_{k} - \beta\varphi (X_{k+1}) + \beta\varphi (X_k) .
\end{array}
\right.
 \end{equation}

Consider the system (\ref{F2}) with arbitrary initial condition $(X_1,Y_1)=(x,y)$.
Observe that if $(x,y)$ is a fixed point for the system (\ref{F2})  then 
\begin{equation}\label{Fp}
(x,y)=\left(a, \frac{\gamma}{4a^2}\right)
 \end{equation}
where 
\[a=\sqrt{\frac{2\beta + \gamma}{2\lambda}},\]
as in (\ref{defa}). 
The linearisation of the system (\ref{F5}) in the neighbourhood of this fixed point is:
\begin{equation}\label{F7}
\left(
\begin{array}{c}
 u_{k+1} \\
v_{k+1} 
\end{array}
\right) = 
\left(
\begin{array}{ll}
 -1 & \left(\varphi^{-1} \left(\frac{y}{\gamma}\right)\right)'\\
  \beta\varphi' (x) &  -1 - \beta \varphi '(x) \left(\varphi^{-1} \left(\frac{y}{\gamma}\right)\right)'
\end{array}
\right) \left(
\begin{array}{c}
 u_{k} \\
v_{k} 
\end{array}
\right),
\end{equation}
where by (\ref{phi4}) and (\ref{Fp})
\begin{equation}\label{F8}
\varphi '(x)=-\frac{2}{x^3}=-\frac{2}{a^3}, \ \ 
\left(\varphi^{-1} \left(\frac{y}{\gamma}\right)\right)'=-\frac{1}{\gamma}\frac{1}{2(y/\gamma)^{3/2}}=-\frac{4a^{3}}{\gamma}.
\end{equation}
Hence, we can rewrite (\ref{F7}) as  
\begin{equation}\label{F9}
\left(
\begin{array}{c}
u_{k+1} \\
v_{k+1} 
\end{array}
\right) = 
\left(
\begin{array}{ll}
-1 & -4\frac{a^{3}}{\gamma}\\
-2\beta a^{-3} &  -1-8\frac{\beta}{\gamma}
\end{array}
\right) \left(
\begin{array}{c}
u_{k} \\
v_{k} 
\end{array}
\right) =: A\left(
\begin{array}{c}
u_{k} \\
v_{k} 
\end{array}
\right).
\end{equation}
Observe, that independent of $a$ (and hence, $\lambda$), $\beta$, and $\gamma\neq 0$, we have for $A$ defined in (\ref{F9})
\begin{equation}\label{detA}
\det A =1.
\end{equation}
The eigenvalues of the matrix $A$ (which do not depend on $\lambda$ neither) 
are 
\begin{equation}\label{F10}
-1<\eta_1 = -1-4\frac{\beta}{\gamma} +\sqrt{16\left(\frac{\beta}{\gamma} \right)^2 +8\frac{\beta}{\gamma}}<0, \ \ \eta_2 = \frac{1}{\eta_1}<-1.
\end{equation}
This implies that if and only if the initial condition $(u_1,v_1)$ for the system (\ref{F9}) is chosen along the eigenvector corresponding to the eigenvalue $\eta_1$, for which $|\eta_1|<1$,
the linearised system (\ref{F7}) converges to the limiting zero state exponentially fast. 

As $a_k$ satisfy (\ref{H2})
and therefore also system (\ref{F1}), our assumption (\ref{As1}) 
is equivalent to that $(X_k=a_k,Y_k)$ does not leave an arbitrarily small neighbourhood of the point (\ref{Fp}), where it is described by the linearised system (\ref{F7}). Hence, by the Hartman-Grobman Theorem exponentially fast convergence of $(u_k,v_k)$, implies
 the same  exponential rate of convergence for the original system $(X_k,Y_k)$ to the limiting state (\ref{Fp}). \hfill$\Box$

\subsection{Gaussian approximation}\label{GL}
First we derive Gaussian approximation for the circular distribution.

\begin{theo}\label{Lexp0}
  For any $N>4$ let  $\bar{X}^{\circ }_{\lambda}$
  be a random vector defined by (\ref{Mt7}) and (\ref{H0c}), and let
  $\bar{a}^{\circ} =(a, \ldots, a)\in {\mathbb{R}}^N$ be the point of minimum of function $H_{\lambda}^{\circ , N}$, i.e.,
	\[
	a=\sqrt{\frac{2\beta + \gamma}{2\lambda}}.
	\]
	
	\noindent
	Assume that $\lambda$ is  of order $N^2$	
	\begin{equation}\label{orderl}
	\lambda  \sim N^2,
	\end{equation}	
	and set
		\begin{equation}\label{N}
	\widetilde{N} = \left[ N^{\frac{1}{7}} \right].
	\end{equation}	
	Then 
  for any continuous  function $g: \mathbb{R}^{\widetilde{N}} \rightarrow \mathbb{R}$ 
with at most polynomial growth, i.e., such that 
	\begin{equation}\label{Ma36}
|g(\bar{x})|\leq P_k(\|\bar{x}\|),\ \ \ \bar{x}\in \mathbb{R} ^{\widetilde{N}},
	\end{equation}
for some polynomial $P_k$ of degree $k\geq 0$ and with constant in $N$ coefficients, one has
	\begin{equation}\label{Ma39}
	{\bf E} g\left(
	{X}^{\circ , N}_{1, \lambda}-{a}, \ldots, {X}^{\circ , N }_{{\widetilde{N}}, \lambda}-{a}
	\right)
	= {\bf E} g\left(
	Z_1, \ldots, Z_{\widetilde{N}}
	\right) 
	\left(1	+
	\ o\left(N^{-\frac{1}{3} }
	\right)
	\right)
	 +o\left( e^{-
		({\log \lambda})^{3/2}}
	\right),
	\end{equation}
where random vector $\left(
Z_1, \ldots, Z_{\widetilde{N}}, \ldots, Z_{2\widetilde{N}}
\right)  $ has a multivariate normal distribution in $\mathbb{R}^{2\widetilde{N}}$ with zero mean vector and covariance
matrix $\Lambda(2\widetilde{N})={\cal H}^{-1}(2\widetilde{N})$, which is inverse of
the
following 
symmetric circulant ${2\widetilde{N}} \times {2\widetilde{N}}$ matrix 
\begin{equation}\label{Hess}
{\cal H}(2\widetilde{N}):= 
\frac{1}{a^3}\left(
\begin{array}{cccccc}
2\beta+\frac{\gamma}{2} & \frac{\gamma}{4} & 0 & \ldots&  0 & \frac{\gamma}{4} \\
\frac{\gamma}{4} & 2\beta+\frac{\gamma}{2}& \frac{\gamma}{4} & 0&  \ldots&  0 \\
\ldots & & & & & \\
\frac{\gamma}{4} & 0 &  \ldots &   0 & \frac{\gamma}{4} & 2\beta+\frac{\gamma}{2}
\end{array}
\right).
\end{equation}
Matrix ${\cal H}(N)$ is the Hessian of function $H_{\lambda}^{\circ, N}$ at point  $(a, \ldots, a)$.
\end{theo}

\bigskip

  The result of this theorem
  gives us a Gaussian  approximation up to a given small error,
  for a vector of ${\widetilde{N}}$ consecutive values
\begin{equation}\label{rv}
  {X}^{\circ, N}_{1+n,  \lambda}, \ldots, {X}^{\circ, N}_{{\widetilde{N}+n}, \lambda}
\end{equation}
  for any $n$ (where we take $(1+n)_{\mod N},  \ldots ,
  ({\widetilde{N}}+n)_{\mod N}$, if needed). Furthermore, $\widetilde{N}\sim N^{1/7}$
  is also  increasing in $N$.

\begin{rem}\label{R2}
  The reduction to a smaller vector (\ref{rv}) of size ${\widetilde{N}}$ is perhaps, only technical. Our approach does not resolve the question whether a similar approximation holds for the entire vector of $N$ entries. The positive answer to this question would certainly simplify further analysis, but would not make substantial changes to the main result on covariance structure.
  \end{rem}

The analogue of Theorem \ref{Lexp0} holds as well for random vector $\bar{X}_{\lambda}$, and we state it below. (However, the inverse matrix  ${\cal H}_1^{-1}$  does not have equally closed form as ${\cal H}^{-1}$ featured in Theorem \ref{Lexp0}.)

\begin{theo}\label{Lexp1}
  Let $\bar{a} =(a_1, \ldots, a_N)\in {\mathbb{R}}^N$ be the point of minimum of function $H_{\lambda}$, and denote ${\cal H}_1(N)$ the Hessian of function $H_{\lambda}$ at $\bar{a}$.
  Under assumptions of Theorem
  \ref{Lexp0} it also holds that 
	\begin{equation}\label{DeMa39}
{\bf E} g\left(
{X}^{ N}_{1, \lambda}-{a}, \ldots, {X}^{ N }_{{\widetilde{N}}, \lambda}-{a}
\right)
= {\bf E} g\left(
Z'_1, \ldots, Z'_{\widetilde{N}}
\right) 
\left(1	+
\ o\left(N^{-\frac{1}{3} }
\right)
\right)
+o\left( e^{-
	({\log \lambda})^{3/2}}
\right),
	\end{equation}
	where random variable $\bar{ \cal Z}'$ has a multivariate normal distribution in $\mathbb{R}^N$ with zero mean vector and covariance matrix ${\cal H}_1^{-1}(N)$.
\end{theo}

\noindent {\bf Proof of Theorem \ref{Lexp0}. }
First we observe that  
\begin{equation}\label{S12}
{\bf E} g\left({X}^{\circ, N}_{1, \lambda}-a, \ldots, {X}^{\circ, N}_{\widetilde{N}, \lambda}-a\right)=\frac{\int_{[0,1]^{{N}}} g(
	x_1-a, \ldots, x_{\widetilde{N}}-a
	)e^{-H_{\lambda}^{\circ, N}(\bar{x})}
	dx_{1} \ldots d x_{{N}}}{
	\int_{[0,1]^{{N}}} 
	e^{-H_{\lambda}^{\circ, N}(\bar{x})}
	dx_{1} \ldots d x_{{N}}
}
\end{equation}

\[=\int_{[0,1]^{\widetilde{N}}} g(x_1-a, \ldots, x_{\widetilde{N}}-a)
f_{{X}^{\circ, N}_{1, \lambda}, \ldots, {X}^{\circ, N}_{\widetilde{N}, \lambda}}(x_{1}, \ldots, x_{\widetilde{N}})dx_{1} \ldots d x_{\widetilde{N}}
\]
 \[=\left(1	+
\ O\left(e^{-\alpha \widetilde{N} }
\right)\right) 
{\bf E} g\left({X}^{\circ, 2\widetilde{N} }_{1, \lambda}-a, \ldots, {X}^{\circ, 2\widetilde{N}}_{\widetilde{N}, \lambda}-a
\right)
,
\]
where the last equality is due to (\ref{S11}).
In turn this yields
\begin{equation}\label{S13}
{\bf E} g\left(
{X}^{\circ, N}_{1, \lambda}-a, \ldots, {X}^{\circ, N}_{\widetilde{N}, \lambda}-a\right)
\end{equation}
\[
=\left(1	+
\ O\left(e^{-\alpha \widetilde{N} }
\right)\right) \frac{\int_{[0,1]^{2\widetilde{N}}} g(x_{1}-a,  \ldots ,x_{\widetilde{N}}-a)e^{-H_{\lambda}^{\circ, 2\widetilde{N}}(
		x_{1},  \ldots ,x_{2\widetilde{N}}
		)}
	dx_{1} \ldots d x_{2\widetilde{N}}}{
	\int_{[0,1]^{2\widetilde{N}}} 
	e^{-H_{\lambda}^{\circ, 2\widetilde{N}}(
		x_{1},  \ldots ,x_{2\widetilde{N}}
		)}
	dx_{1} \ldots d x_{2\widetilde{N}}
}.
\]

For simplicity of notation we begin with considering $\bar{X}^{\circ}_{\lambda}=\bar{X}^{\circ, N}_{\lambda}$, and then  turn to $\bar{X}^{\circ, 2\widetilde{N} }_{\lambda}$.

Whenever it does not cause a confusion, for any $N'\geq \widetilde{N}$  
	and $\bar{x}=(x_1, \ldots, x_{N'})$ we shall write
	\[g(\bar{x})= g(x_1, \ldots, x_{\widetilde{N}}).\]

For any $\bar{x}=(x_1, \ldots, x_N)$  set 
$x_{N+1}=x_1$, and rewrite (\ref{H0c}) as
\[H_{\lambda}^{\circ , N}(\bar{x})= \sum_{k=1}^{N} \frac{\beta}{x_k} + \sum_{k=1}^{N}\frac{\gamma}{x_k+x_{k+1}} +\lambda \sum_{k=1}^{N} x_k.\]
For the point of minimum of $H_{\lambda}^{\circ , N}(\bar{x})$ defined by (\ref{defa}) under the assumption that $\lambda=\lambda(N) \rightarrow \infty$, it holds that
\[\bar{a}^{\circ}=(a, \ldots, a) \in (0,1)^N. \]
Observe that $a=a(\lambda)$ does not depend on $N$ (except via $\lambda$).
Consider the Taylor expansion of $H_{\lambda}^{\circ , N}(\bar{x}) $ in the neighbourhood 
 \begin{equation}\label{B}
 B_{\frac{\varepsilon}{\sqrt{\lambda}}}(\bar{a}^{\circ}):=\left\{ \bar{x} : \max_{1\leq i\leq N}
|x_i-a|
< \frac{\varepsilon}{\sqrt{\lambda}} \right\}
\end{equation}
of  $\bar{a}^{\circ}$. For all  $\bar{x} \in \left[ B_{\frac{\varepsilon}{\sqrt{\lambda}}}(\bar{a}^{\circ})\right]=\left\{ \bar{x} : \max_{1< i\leq N}
|x_i-a|
\leq \frac{\varepsilon}{\sqrt{\lambda}} \right\}$ we have
\begin{equation}\label{An2}
H_{\lambda}^{\circ , N}(\bar{x})= H_{\lambda}^{\circ , N}(\bar{a}^{\circ}) +
\frac{1}{2}\sum_{k=1}^{N}\sum_{l=1}^{N}\frac{\partial ^2H_{\lambda}^{\circ}}{\partial x_k\partial x_l}(\bar{a}^{\circ})
  (x_k-a)(x_l-a) 
\end{equation}
\[+ 
\frac{1}{3!}\sum_{k=1}^{N}\sum_{l=1}^{N}\sum_{j=1}^{N}\frac{\partial ^3H_{\lambda}^{\circ}}{\partial x_k\partial x_l\partial x_j}(\bar{y}) \,
  (x_k-a)(x_l-a)(x_j-a),
\]
where 
\[ \bar{y}= \bar{y}(\bar{x}) \in B_{\frac{\varepsilon}{\sqrt{\lambda}}}(\bar{a}^{\circ}).
\]
Compute now
\begin{equation}\label{An3}
\frac{\partial ^2H_{\lambda}^{\circ}}{\partial x_k\partial x_l}(\bar{a}^{\circ})= 
\left\{
\begin{array}{ll}
(2\beta+\frac{\gamma}{2})a^{-3}, & \ k=l, \\ \\
\frac{\gamma}{4a^3}, & \ |k-l|=1, \\ \\
0, & \ |k-l|>1. \\
\end{array}
\right.
\end{equation}
Hence, the  Hessian of function  $H_{\lambda}^{\circ , N}(\bar{x}) $  at the point $\bar{a}^{\circ}$ is indeed matrix ${\cal H}={\cal H}(N)$ defined in (\ref{Hess}) as claimed in the theorem.

The eigenvalues of ${\cal H}(N)$, call them  $\nu^{\circ}_j, j=1, \ldots, N,$   are given by 
\begin{equation}\label{eigenM}
\nu^{\circ}_j=\frac{1}{a^3}\left(\left( 2\beta+\frac{\gamma}{2} \right) +2\frac{\gamma}{4} \cos \left( 2\pi \frac{ j}{N}\right)\right)  
\end{equation}
(see, e.g., \cite{D}). The latter 
yields a uniform bound
\begin{equation}\label{eigen}
 \frac{2\beta}{a^3} \leq \nu^{\circ}_j \leq  \frac{2\beta+\frac{\gamma}{2}}{a^3}\ \ \ \  j=1, \ldots, N.
\end{equation}
which  proves that ${\cal H}$ is positive-definite (it confirms in particular,  that point $\bar{a}^{\circ}$ is the minimum). 

For a later reference let us point out a useful corollary of the assertion (\ref{eigen}): for any $\bar{x}\in \mathbb{R}^N$
\begin{equation}\label{n6}
\frac{2\beta}{a^3} \,  \|\bar{x}\|²\leq  \bar{x}'{\cal H}\bar{x}\leq \frac{2\beta+\frac{\gamma}{2}}{a^3} \, \|\bar{x}\|².
\end{equation}

It is convenient to use here the vector form. We denote
$\bar{x}$ the {\it row}-vector, and $\bar{x}'$  the {\it column}-vector. 
Then we get from (\ref{An2}) taking into account the third-order derivatives as well
\begin{equation}\label{An4}
H_{\lambda}^{\circ , N}(\bar{x})= H_{\lambda}^{\circ , N}(\bar{a}^{\circ}) +
\left( 1+O\left(\frac{\varepsilon}{a \sqrt{\lambda}}\right) \right) \frac{1}{2}(\bar{x}-\bar{a}^{\circ}){\cal H}(\bar{x}-\bar{a}^{\circ})'.
\end{equation}

Let function $g(\bar{x}): \mathbb{R}^N \rightarrow \mathbb{R}$ satisfy the assumptions of Theorem \ref{Lexp0}. To make use of   (\ref{An4}) in computing integral   
\begin{equation}\label{IA}
\int_{[0,1]^N }  g(\bar{x}-\bar{a}^{\circ})
 e^{- H_{\lambda}^{\circ , N}(\bar{x})  }d\bar{x},
\end{equation}
featured in (\ref{S12}),
we shall split its area of integration $[0,1]^N $ as follows.
First we fix  $0<\varepsilon <1$ arbitrarily and define for any $\bar{x}\in [0,1]^N$
the set of its components  which deviate from $a$ at least by $\frac{\varepsilon}{\sqrt{\lambda}}$:
\begin{equation}\label{a11}
{\cal K}(\bar{x}) = \{i: |x_i-a|\geq \frac{\varepsilon}{\sqrt{\lambda}} \}.
\end{equation}
Then for any ${\cal K} \subseteq \{1, \ldots, N\} $ define a subset of vectors $\bar{x}\in [0,1]^N$ 
for which ${\cal K}(\bar{x})=K$
 \begin{equation}\label{au1}
 {S}_{\cal K} = \{\bar{x}\in [0,1]^N: {\cal K}(\bar{x}) = {\cal K} \}.
 \end{equation}
 In particular,  in this notation we have
 \[B_{
 	\frac{\varepsilon}{\sqrt{\lambda}}
 (\bar{a}^{\circ})
 } ={S}_{\emptyset}.\]
This gives us a decomposition
\[
[0,1]^N =\cup_{{\cal K}\subseteq \{1, \ldots, N\}} {S}_{\cal K},
	\]
which allows us to split
\[
\int_{[0,1]^N }  
g(\bar{x}-\bar{a}^{\circ}) e^{- H_{\lambda}^{\circ , N}(\bar{x})  }d\bar{x}
\]
\begin{equation}\label{a16}
= \sum_{{\cal K}\subseteq \{1, \ldots, N\}}\int_{ \{ \bar{x}\in [0,1]^N: \ {\cal K}(\bar{x}) ={\cal K}\}}  
g(\bar{x}-\bar{a}^{\circ}) e^{- H_{\lambda}^{\circ , N}(\bar{x})  }d\bar{x}.	
\end{equation}

Given a non-empty set ${\cal K}\subset \{1, \ldots, N\}$ let ${\cal H}_{{\cal K}}$ denote a matrix obtained from ${\cal H}$ by deleting  all the columns and rows with indices in ${\cal K}$. Alternatively, ${\cal H}_{{\cal K}}$ can be defined as a matrix of the quadratic form obtained from $\bar{x}{\cal H}\bar{x}' $ by deleting all the terms which include $x_i$ with $i\in {\cal K}$. (When convenient we shall also write ${\cal H}_{{\emptyset}}={\cal H}.$)
Let $\bar{x}_{{\cal K}}$ denote  an  $(N-|{\cal K}|)$-dimensional vector obtained from $\bar{x}$ by deleting the components with indexes in ${\cal K}=\{i_1, \ldots, i_{|{\cal K}|}\}$, where $i_1<\ldots< i_{|{\cal K}|}$.
Define a set of differences between the consecutive $i_k$ in this set:
\[m_k= i_{k+1}-i_k-1, k=1, \ldots,  |{\cal K}|-1,  \]
\[m_{|{\cal K}|}= N-i_{|{\cal K}|}+i_1-1.\]
Notice that by this definition
\[\sum_{1\leq k \leq|{\cal K}|: m_k>0}m_k= N-|{\cal K}|.\]
 Then 
$\bar{x}_{{\cal K}}{\cal H}_{{\cal K}}\bar{x}'_{{\cal K}} $ can be written as a sum of
independent quadratic forms, each corresponding to  $m_k\times m_k$ matrix $A_{m_k}$  defined by
	\begin{equation}\label{HessK} 
A_{m}:=	\frac{1}{a^3}\left(
	\begin{array}{cccccc}
          2\beta+\frac{\gamma}{2} &
                                    \frac{\gamma}{4} &
                                                       0 &
                                                           \ldots &
                                                                    0 &
                                                                        0 \\
	\frac{\gamma}{4} & 2\beta+\frac{\gamma}{2}& \frac{\gamma}{4} & 0&  \ldots&  0 \\
	\ldots & & & & & \\
	0 & 0 &  \ldots &   0 & \frac{\gamma}{4} & 2\beta+\frac{\gamma}{2}
	\end{array}
	\right),
      \end{equation}
      when $m>1$, and also for $m=1$ set
      \[A_1:=	\frac{1}{a^3}\left( 2\beta+\frac{\gamma}{2} \right). \]
Hence, 
\begin{equation}\label{detHk}
\det {\cal H}_{{\cal K}} = \prod_{1\leq k \leq|{\cal K}|: m_k>0}\det A_{m_k}.
\end{equation}
      The eigenvalues of 3-diagonal Toeplitz matrix $A_m$ are well-known (see, e.g., \cite{K}), these are
\begin{equation}\label{eiv}
  c_j(m)=\frac{1}{a^3} \left(
    \left( 2\beta+\frac{\gamma}{2}\right)+
    2\frac{\gamma}{4}\cos \left( \pi \frac{ j}{m+1}\right)
  \right), \ \ \ j=1, \ldots, m,
\end{equation}
yielding uniform bounds as in (\ref{eigen})
\begin{equation}\label{eivb}
  \frac{2\beta}{a^3} \leq c_j(m) \leq
  \frac{1}{a^3} \left( 2\beta+{\gamma} \right), \ \ \ j=1, \ldots, m.
\end{equation}
Furthermore, this implies that the eigenvalues of matrix ${\cal H}_{{\cal K}}$
as well are  between same values $\frac{2\beta}{a^3}$   and $ \frac{1}{a^3} \left( 2\beta+{\gamma} \right)$, and we have as in  (\ref{n6})
\begin{equation}\label{n6*}
\frac{2\beta}{a^3} \,  \|\bar{x}_{{\cal K}}\|²\leq  \bar{x}_{{\cal K}}'{\cal H}_{{\cal K}}\bar{x}_{{\cal K}}\leq \frac{2\beta+{\gamma}}{a^3} \, \|\bar{x}_{{\cal K}}\|²
\end{equation}
for any $\bar{x}_{{\cal K}}\in \mathbb{R}^{N-|{{\cal K}}|}$.

This together with (\ref{detHk}) confirms that  matrix ${\cal H}_{{\cal K}}$ is also positive-definite for all ${\cal K}\subset \{1, \ldots, N\}$ (if fact, it is positive-definite since ${\cal H}$ is). Note here for further reference that the bounds (\ref{eivb}) together with (\ref{detHk}) imply
\begin{equation}\label{detHkb}
\left( \frac{2\beta}{a^3} \right)^{N-|{\cal K}|}\leq \det {\cal H}_{{\cal K}} \leq \left( \frac{1}{a^3} \left( 2\beta+{\gamma} \right) \right)^{N-|{\cal K}|}
\end{equation}
for any set ${\cal K}$.

Fix a non-empty set ${\cal K}\subset \{1, \ldots, N\}$ arbitrarily.
Note that for any $\bar{x} \in \left[ B_{\frac{\varepsilon}{\sqrt{\lambda}}}(\bar{a}^{\circ})\right]$, where  the  approximation (\ref{An4}) is still valid,  and such that 
${\cal K}(\bar{x}) ={\cal K}$, i.e., 
\[
\{i: |x_i-a|= \frac{\varepsilon}{\sqrt{\lambda}} \}={\cal K},
\]
we have for any $i\in {\cal K}$ 
\[\left(2\beta+\frac{\gamma}{2} \right)(x_i-a_i)^2+2\frac{\gamma}{4} (x_{i+1}-a_{i+1})(x_{i}-a_i)+2\frac{\gamma}{4} (x_{i-1}-a_{i-1})(x_{i}-a_i)\]
\[\geq \left(\frac{\varepsilon}{\sqrt{\lambda}} \right)^2\left(2\beta+\frac{\gamma}{2} -\gamma \right), \]
where we replace $i\pm 1$ by $(i\pm 1)_{\mbox{\small{mod }} N}$ if needed.
Therefore making use of the last bound we derive for all such $\bar{x} $  and 
 small $\varepsilon$
\[
H_{\lambda}^{\circ , N}(\bar{x})=H_{\lambda}^{\circ , N}(\bar{a}^{\circ}) +
\left( 1+O\left(\frac{\varepsilon}{a \sqrt{\lambda}}\right) \right) \frac{1}{2}(\bar{x}-\bar{a}^{\circ}){\cal H}(\bar{x}-\bar{a}^{\circ})'\]
\begin{equation}\label{Ja15}
 \geq H_{\lambda}^{\circ , N}(\bar{a}^{\circ})+\left( 1+O\left(\frac{\varepsilon}{a \sqrt{\lambda}}\right) \right) \frac{1}{2}
 |{\cal K}| \frac{1}{a^3}\left( \frac{\varepsilon}{\sqrt{\lambda}}\right) ^2\left(2\beta-\frac{\gamma}{2} \right) 
 \end{equation}
 \[+  \left( 1+O\left(\frac{\varepsilon}{a \sqrt{\lambda}}\right) \right) \frac{1}{2}(\bar{x}_{{\cal K}} -\bar{a}^{\circ}_{{\cal K}}) {\cal H}_{{\cal K}}(\bar{x}_{{\cal K}} -\bar{a}^{\circ}_{{\cal K}})'
.\]

By the convexity of function $H_{\lambda}^{\circ, N}$ the last bound holds as well for all 
$\bar{x} \in [0,1]^N$  with ${\cal K}(\bar{x})={\cal K}$, 
and therefore for a continuous function $g$  we get 
\begin{equation}\label{a12}
\left| \int_{
	\{ 
	\bar{x} \in [0,1]^N : \ {\cal K}(\bar{x})={\cal K}
	\}
}
g(\bar{x}-\bar{a}^{\circ})	e^{-H_{\lambda}^{\circ, N}(\bar{x})   }d\bar{x}\right| 
\end{equation}
\[
\leq  Ce^{-H_{\lambda}^{\circ, N}(\bar{a}^{\circ})-|{\cal K}| B {\varepsilon}^2{\sqrt{\lambda}}
}
\
 \int_{
	\bar{x}_{{\cal K}}: \max_i|x_i-a|\leq \frac{\varepsilon}{\sqrt{\lambda}}
}	e^{-
\left( 1+O\left(\frac{\varepsilon}{a \sqrt{\lambda}}\right) \right) \frac{1}{2}(\bar{x}_{{\cal K}} -\bar{a}^{\circ}_{{\cal K}}) {\cal H}_{{\cal K}}(\bar{x}_{{\cal K}} -\bar{a}^{\circ}_{{\cal K}})'
}d\bar{x}_{{\cal K}},\]
where $$C=\max_{\bar{x} \in [0,1]^N } |g(\bar{x}-\bar{a}^{\circ})|, \ \ B= \frac{2\beta-\frac{\gamma}{2}}{3\left(\beta + \frac{\gamma}{2}\right)^{3/2}}.$$
This yields  a straightforward bound
\begin{equation}\label{Ja12}
\left| \int_{
	\{ 
	\bar{x} \in [0,1]^N : \ {\cal K}(\bar{x})={\cal K}
	\}
}
g(\bar{x}-\bar{a}^{\circ})	e^{-H_{\lambda}^{\circ, N}(\bar{x})   }d\bar{x}\right| 
\end{equation}
\[
\leq  Ce^{-H_{\lambda}^{\circ, N}(\bar{a}^{\circ})-|{\cal K}| B {\varepsilon}^2{\sqrt{\lambda}}
	}
	\
	\frac{\left( (2\pi\left( 1+O\left(\frac{\varepsilon}{a \sqrt{\lambda}}\right) \right)\right)^{\frac{N-|{\cal K}|}{2}}}{\sqrt{\det {\cal H}_{{\cal K}}}}
.\]

Making use of (\ref{Ja12}) we derive
\[
\left|\int_{[0,1]^N \setminus B_{\frac{\varepsilon}{\sqrt{\lambda}}}(\bar{a}^{\circ})}  g(\bar{x}-\bar{a}^{\circ})
  e^{- H_{\lambda}^{\circ, N}(\bar{x})  }d\bar{x}\right|
\leq \int_ {\left\{ \bar{x}\in [0,1]^N:
\min_{1\leq i\leq N}|x_i-a|
\geq  \frac{\varepsilon}{\sqrt{\lambda}} \right\}}
  C
  e^{- H_{\lambda}^{\circ, N}(\bar{x})  }d\bar{x}
\]
\begin{equation}\label{Ja17}
	+ C\sum_{K=1}^{N-1} \ \sum_{{\cal K}\subseteq \{1, \ldots, N\}: |{\cal K}|
	=K}\int_{ \{ \bar{x}\in [0,1]^N \setminus B_{\frac{\varepsilon}{\sqrt{\lambda}}}(\bar{a}^{\circ}): \ {\cal K}(\bar{x}) ={\cal K}\}}  
e^{- H_{\lambda}^{\circ , N}(\bar{x})  }d\bar{x}
	\end{equation}

\[
\leq C
e^{-H_{\lambda}^{\circ , N}(\bar{a}^{\circ})}
e^{-N B {\varepsilon}^2{\sqrt{\lambda}}}
 \]
  \[ +C
 e^{-H_{\lambda}^{\circ ,N}(\bar{a}^{\circ})
 }
 \
 \frac{ 
 	(2\pi)^{ \frac{N}{2} } 
 }{
 	\sqrt{\det {\cal H}} 
 }
 \
 \sum_{K=1}^{N-1} \ \sum_{{\cal K}\subseteq \{1, \ldots, N\}: |{\cal K}|
 	=K}
 e^{-K B {\varepsilon}^2{\sqrt{\lambda}}
 }
 \frac{
 	\sqrt{\det 
 		{\cal H}
 }}{\sqrt{\det 
 		{\cal H}_{\cal K}
 }}
\left( 1+O\left(\frac{\varepsilon}{a \sqrt{\lambda}}\right) \right)^{\frac{N-|{\cal K}|}{2}}
.
 \]

 \bigskip

 \noindent
 {\bf Claim.}
For any set ${\cal K}\subseteq \{1, \ldots, N\}$
\begin{equation}\label{Cl1}
\frac{
	\sqrt{\det 
		{\cal H}
	}}{
\sqrt{\det 
{\cal H}_{\cal K}
}}\leq \sqrt{2}\left(\frac{1}{a^3}\left(2\beta+\frac{\gamma}{2} \right)\right)^{|{\cal K}|/2}.
\end{equation}        
 
\noindent
 {\bf Proof of the Claim.} It follows from the definition of circular matrix ${\cal H}$ (see (\ref{Hess})) that for any $i \in \{1, \ldots, N\}$
\begin{equation}\label{De1}
\det 
{\cal H}\leq \frac{1}{a^3}\left(2\beta+\frac{\gamma}{2} \right)\det 
{\cal H}_{\{i\}}+2\left(\frac{\gamma}{4a^3} \right)^N,
\end{equation}
where
\[\det 
{\cal H}_{\{i\}}= \det A_{N-1},\]
as defined in (\ref{HessK}).
Then for $A_m$ we also derive from the definition (\ref{HessK})
\begin{equation}\label{De3}
\det A_m
\leq \frac{1}{a^3}\left(2\beta+\frac{\gamma}{2} \right)\det A_l \det A_{m-l-1},
\end{equation}
for any $1\leq l<m$. When $m=N-1$ this is equivalent (see also (\ref{detHk})) to 
\begin{equation}\label{De2}
\det 
{\cal H}_{\{i\}}\leq \frac{1}{a^3}\left(2\beta+\frac{\gamma}{2} \right)\det 
{\cal H}_{\{i,j\}}
\end{equation}
for any $j\neq i$. Hence, by the recurrence we get
for any set ${\cal K}\subseteq \{1, \ldots, N\}$ 
\begin{equation}\label{De4}
\det 
{\cal H}_{\{i\}}\leq \left(\frac{1}{a^3}\left(2\beta+\frac{\gamma}{2} \right)\right)^{|{\cal K}|-1}\det 
{\cal H}_{{\cal K}}. 
\end{equation}
Making use of the last bound in (\ref{De1}) we obtain
\begin{equation}\label{De5}
\det 
{\cal H}\leq \left(\frac{1}{a^3}\left(2\beta+\frac{\gamma}{2} \right)\right)^{|{\cal K}|}\det 
{\cal H}_{{\cal K}}+2\left(\frac{\gamma}{4 a^3 } \right)^N.
\end{equation}
Then taking into account the assumption that $\gamma \leq \beta$, 
we get using the lower bound in (\ref{detHkb})
that 
\[2\left(\frac{\gamma}{4 a^3 } \right)^N \leq \left(\frac{1}{a^3}\left(2\beta+\frac{\gamma}{2} \right)\right)^{|{\cal K}|}
\left(\frac{2\beta}{a^3}\right)^{N-|{\cal K}|}
\leq \left(\frac{1}{a^3}\left(2\beta+\frac{\gamma}{2} \right)\right)^{|{\cal K}|}\det 
{\cal H}_{{\cal K}}.\]
Combining the last bound with (\ref{De5}), we get
\begin{equation}\label{De6}
\det 
{\cal H}\leq 2\left(\frac{1}{a^3}\left(2\beta+\frac{\gamma}{2} \right)\right)^{|{\cal K}|}\det 
{\cal H}_{{\cal K}}
\end{equation}
for any set ${\cal K}\subseteq \{1, \ldots, N\}$, which proves the Claim. \hfill$\Box$
\bigskip

The derived  bound (\ref{Cl1})  allows us to obtain from   (\ref{Ja17}):
          \begin{equation}\label{Ja18}
\left|\int_{[0,1]^N \setminus B_{\frac{\varepsilon}{\sqrt{\lambda}}}(\bar{a}^{\circ})}  g(\bar{x}-\bar{a}^{\circ})
  e^{- H_{\lambda}^{\circ , N}(\bar{x})  }d\bar{x}\right|
\end{equation}
\[
\leq C
e^{-H_{\lambda}^{\circ , N}(\bar{a}^{\circ})}
e^{-N B {\varepsilon}^2{\sqrt{\lambda}}}
\]
\[ +C
e^{-H_{\lambda}^{\circ , N}(\bar{a}^{\circ})
}
\
\frac{ 
	(2\pi)^{ \frac{N}{2} } 
}{
	\sqrt{\det {\cal H}} 
}
\
\sum_{K=1}^{N-1} 
\left(
\begin{array}{c}
N \\
K
\end{array}
\right)
e^{-K B {\varepsilon}^2{\sqrt{\lambda}}
}
\sqrt{2}\left(\frac{1}{a^3}\left(2\beta+\frac{\gamma}{2} \right)\right)^{|{\cal K}|/2}\left( 1+O\left(\frac{\varepsilon}{a \sqrt{\lambda}}\right) \right)^{\frac{N-|{\cal K}|}{2}}
\]
\[
\leq C
e^{-H_{\lambda}^{\circ , N}(\bar{a}^{\circ})}
e^{-N B {\varepsilon}^2{\sqrt{\lambda}}}
\]
\[ +2C
e^{-H_{\lambda}^{\circ , N}(\bar{a}^{\circ})
}
\
\frac{ 
	(2\pi)^{ \frac{N}{2} } 
}{
	\sqrt{\det {\cal H}} 
}
\left(
\left(
e^{- B {\varepsilon}^2{\sqrt{\lambda}}
}
\frac{\sqrt{2\beta+\frac{\gamma}{2} }}{a^{3/2}}
+1
\right)^{N}
-1\right)\left( 1+O\left(\varepsilon\right) \right)^{\frac{N}{2}},
\]
where we took into account that $a\sqrt{\lambda} \sim const$.
Let us set from now on
\begin{equation}\label{Ja20}
\varepsilon = 	\frac{\log \lambda}{\lambda^{\frac{1}{4}}}.
\end{equation}
This choice together with the upper bound in  (\ref{detHkb}) yields 
for all 
\begin{equation}\label{Ja20l}
\lambda >N
\end{equation}
that 
\[e^{-N B {\varepsilon}^2{\sqrt{\lambda}}}\frac
{
	\sqrt{\det {\cal H}} 
}{ 
	(2\pi)^{ \frac{N}{2} } 
}\leq \left(e^{- B ({\log \lambda})^2} \frac{2\beta +\gamma}{a^3\sqrt{2\pi}}\right)^N= O\left(e^{- \frac{B}{2} ({\log \lambda})^2}  \right)^N,
\]
which gives us a bound
\[\left(
e^{- B {\varepsilon}^2{\sqrt{\lambda}}
}
\frac{\sqrt{2\beta+\frac{\gamma}{2} }}{a^{3/2}}
+1
\right)^{N}
-1 \leq  O\left(Ne^{- \frac{B}{2} ({\log \lambda})^2}  \right)= O\left(e^{- \frac{B}{3} ({\log \lambda})^2}  \right). \]
Substituting the last two bounds into
(\ref{Ja18}) allows us to derive for all $\lambda >N$
\begin{equation}\label{Ja21}
\left| \int_{[0,1]^N \setminus B_{\frac{\varepsilon}{\sqrt{\lambda}}}(\bar{a}^{\circ})}  
g(\bar{x}-\bar{a}^{\circ})e^{- H_{\lambda}^{\circ, N}(\bar{x})  }d\bar{x}	\right|=e^{-H_{\lambda}^{\circ, N}(\bar{a}^{\circ})
}
\
\frac{ 
	(2\pi)^{ \frac{N}{2} } 
}{
	\sqrt{\det {\cal H}} 
}O\left(e^{- \frac{B}{4} ({\log \lambda})^2}  \right)
\left( 1+O\left(\varepsilon\right) \right)^{\frac{N}{2}}
.
	\end{equation}

	Repeating the same arguments as in (\ref{Ja15})-(\ref{Ja21}), and keeping choice (\ref{Ja20}) we get a similar to (\ref{Ja21}) bound when $H_{\lambda}^{\circ, N}$ is replaced by the corresponding quadratic form ${\cal H}$, namely 
	\begin{equation}\label{Ja22}
	\int_{B_{1}(\bar{a}^{\circ}) \setminus B_{\frac{\varepsilon}{2\sqrt{\lambda}}}(\bar{a}^{\circ})}  g(\bar{x}-\bar{a}^{\circ})
	e^{- \frac{1}{2}(\bar{x}-\bar{a}^{\circ}){\cal H}(\bar{x}-\bar{a}^{\circ})'  }d\bar{x}= 	
	\frac{ 
		(2\pi)^{ \frac{N}{2} } 
	}{
		\sqrt{\det {\cal H}} 
	}
	\ o\left(e^{-\frac{B}{4} ({\log \lambda})^2
	}\right)\left( 1+O\left(\varepsilon\right) \right)^{\frac{N}{2}}.
    \end{equation}
    Notice that in the last integral a ball $B_{1}(\bar{a}^{\circ})$ replaces a box
    $[0,1]^N$ in (\ref{Ja21}), but the former arguments are still applicable.

Next  for any 
$\bar{x}\in B_{\frac{ \varepsilon}{\sqrt{\lambda}}}(\bar{a}^{\circ}) \cap \{\bar{x}: \|\bar{x}-\bar{a}^{\circ}\|\leq \varepsilon \lambda^{-\frac{1}{2}+q}\}$, where $0<q<\frac{1}{8}$, we derive with a help of (\ref{n6}) 
\begin{equation}\label{De9}
O(\varepsilon)
(\bar{x}-\bar{a}^{\circ}){\cal H}(\bar{x}-\bar{a}^{\circ})'  
\leq O\left(\varepsilon a^{-3}\|\bar{x}-\bar{a}^{\circ}\|^2\right)=O\left(\varepsilon^3 a^{-3}\frac{1}{\lambda^{1-2q}}\right)
\end{equation}
\[= O\left(\lambda^{-\frac{1}{4}+2q}(\log \lambda)^3\right).\]
This gives us the following approximation
\begin{equation}\label{De11}
\int_{B_{\frac{ \varepsilon}{\sqrt{\lambda}}}(\bar{a}^{\circ}) \cap \{\bar{x}: \|\bar{x}-\bar{a}^{\circ}\|\leq \varepsilon \lambda^{-\frac{1}{2}+q}\}
}  
g\left({\bar{x}-\bar{a}^{\circ}}\right)e^{- (1+O(\varepsilon))\frac{1}{2}(\bar{x}-\bar{a}^{\circ}){\cal H}(\bar{x}-\bar{a}^{\circ})'    }d\bar{x}
\end{equation}
\[=\left(1 + O\left(\lambda^{-\frac{1}{4}+2q}(\log \lambda)^3\right)
\right)
\int_{B_{\frac{ \varepsilon}{\sqrt{\lambda}}}(\bar{a}^{\circ}) \cap \left\{\bar{x}: \|\bar{x}-\bar{a}^{\circ}\|\leq \varepsilon \lambda^{-\frac{1}{2}+q}\right\}
}  
g\left({\bar{x}-\bar{a}^{\circ}}\right)e^{- \frac{1}{2}(\bar{x}-\bar{a}^{\circ}){\cal H}(\bar{x}-\bar{a}^{\circ})'    }d\bar{x}.\]
On the other hand, again due to (\ref{n6}), we have for some positive constants $c,c_1$
\begin{equation}\label{De8}
\int_{B_{\frac{ \varepsilon}{\sqrt{\lambda}}}(\bar{a}^{\circ}) \cap \left\{\bar{x}: \|\bar{x}-\bar{a}^{\circ}\|\geq \varepsilon \lambda^{-\frac{1}{2}+q}\right\}
}  
g\left({\bar{x}-\bar{a}^{\circ}}\right)e^{- (1+O(\varepsilon))\frac{1}{2}(\bar{x}-\bar{a}^{\circ}){\cal H}(\bar{x}-\bar{a}^{\circ})'    }d\bar{x}
\end{equation}
\[
\leq C	
\int_{
	B_{\frac{ \varepsilon}{\sqrt{\lambda}}}(\bar{a}^{\circ}) \cap \left\{\bar{x}: \|\bar{x}-\bar{a}^{\circ}\|\geq \varepsilon \lambda^{-\frac{1}{2}+q}\right\}
}  
e^{- a^{-3}c\|\bar{x}-\bar{a}^{\circ}\|^2}d\bar{x} 
\leq O\left(\lambda^{-\frac{3}{4}}\right)^N
e^{-c_1\lambda^{3/2}\varepsilon^2
	\lambda^{-1+2q}}
\]
\[
=  
\frac{(2\pi)^{\frac{N}{2}}}{\sqrt{\det{\cal H}}}
o\left( 
e^{-c_1
	\lambda^{q}}
\right).
\]
Combining (\ref{De8}) and (\ref{De11}) we obtain the following bound
\begin{equation}\label{De12}
\int_{ B_{\frac{\varepsilon}{\sqrt{\lambda}}}(\bar{a}^{\circ})}  
g(\bar{x}-\bar{a}^{\circ})e^{- H_{\lambda}^{\circ , N}(\bar{x})  }d\bar{x}=e^{-H_{\lambda}^{\circ , N}(\bar{a}^{\circ})
}\int_{
	B_{\frac{ \varepsilon}{\sqrt{\lambda}}}(\bar{a}^{\circ}) 
}  
g\left({\bar{x}-\bar{a}^{\circ}}\right)e^{- (1+O(\varepsilon))\frac{1}{2}(\bar{x}-\bar{a}^{\circ}){\cal H}(\bar{x}-\bar{a}^{\circ})'    }
d\bar{x}
\end{equation}
\[=e^{-H_{\lambda}^{\circ , N}(\bar{a}^{\circ})
}
\int_{B_{\frac{ \varepsilon}{\sqrt{\lambda}}}(\bar{a}^{\circ}) \cap \left\{\bar{x}: \|\bar{x}-\bar{a}^{\circ}\|\leq \varepsilon \lambda^{-\frac{1}{2}+q}\right\}
}  
g\left({\bar{x}-\bar{a}^{\circ}}\right)e^{- \frac{1}{2}(\bar{x}-\bar{a}^{\circ}){\cal H}(\bar{x}-\bar{a}^{\circ})'    }d\bar{x}\left(1+
o\left( \lambda^{-\frac{1}{4}+3q}
\right)\right)\]
\[
+\frac{(2\pi)^{\frac{N}{2}}e^{-H_{\lambda}^{\circ , N}(\bar{a}^{\circ})
}}{\sqrt{\det{\cal H}}}o\left( 
e^{-c_1
	\lambda^{q}}
\right)
\]
\[=e^{-H_{\lambda}^{\circ , N}(\bar{a}^{\circ})
}
\int_{B_{\frac{ \varepsilon}{\sqrt{\lambda}}}(\bar{0}) 
}  
g\left({\bar{x}}\right)e^{- \frac{1}{2}\bar{x}
	{\cal H}\bar{x}'    }d\bar{x}
\left(1+
o\left( \lambda^{-\frac{1}{4}+3q}
\right)\right)+\frac{(2\pi)^{\frac{N}{2}}e^{-H_{\lambda}^{\circ , N}(\bar{a}^{\circ})
}}{\sqrt{\det{\cal H}}}o\left( 
e^{-c_1
	\lambda^{q}}
\right)
\]
for arbitrary $0<q<\frac{1}{8}$. This together with (\ref{Ja21}) yields
\begin{equation}\label{23Ja1}
\int_{[0,1]^N }  
g(\bar{x}-\bar{a}^{\circ})e^{- H_{\lambda}^{\circ , N}(\bar{x})  }d\bar{x}	
\end{equation}
\[
=e^{-H_{\lambda}^{\circ, N}(\bar{a}^{\circ})
}
\int_{B_{\frac{ \varepsilon}{\sqrt{\lambda}}}(\bar{0}) 
}  
g\left({\bar{x}}\right)e^{- \frac{1}{2}\bar{x}{\cal H}\bar{x}'    }d\bar{x}\left( 1
+
o\left( \lambda^{-\frac{1}{4}+3q}
\right) \right)\]
\[
+
e^{-H_{\lambda}^{\circ, N}(\bar{a}^{\circ})
}\frac{(2\pi)^{\frac{N}{2}}}{\sqrt{\det{\cal H}}}
o\left( e^{-\frac{B}{3}
	({\log \lambda})^2}
\right)\left( 1+O\left(\varepsilon\right) \right)^{\frac{N}{2}}
\]
for arbitrary  $0<q<\frac{1}{8}$.  Let us rewrite the latter as 
\begin{equation}\label{Ja25}
\int_{[0,1]^N }  
g(\bar{x}-\bar{a}^{\circ})e^{- H_{\lambda}^{\circ, N}(\bar{x})  }d\bar{x}	
\end{equation}
\[
=e^{-H_{\lambda}^{\circ, N}(\bar{a}^{\circ})
}
\left( 1
+
o\left( \lambda^{-\frac{1}{4}+3q}
\right)
\right)
\]
\[\times\left( 
\int_{{\bf R}^N
}  
g(\bar{x}) e^{- \frac{1}{2}\bar{x}{\cal H}\bar{x}'    }d\bar{x} 
-
\int_{{\bf R}^N \setminus  B_{1}(\bar{0})} g(\bar{x})
e^{- \frac{1}{2}\bar{x}{\cal H}\bar{x}' }d\bar{x}+
\int_{B_{1}(\bar{0})
	\setminus B_{\frac{ \varepsilon}{\sqrt{\lambda}}}(\bar{0})
} g(\bar{x})
e^{- \frac{1}{2}\bar{x}{\cal H}\bar{x}' }d\bar{x}
\right)
\]
\[+
e^{-H_{\lambda}^{\circ, N}(\bar{a}^{\circ})
}\frac{(2\pi)^{\frac{N}{2}}}{\sqrt{\det{\cal H}}}
o\left( e^{-\frac{B}{4}
	({\log \lambda})^2}
\right)\left( 1+O\left(\varepsilon\right) \right)^{\frac{N}{2}}
.
\]
By the property (\ref{eigen}) and the assumption on $g$ we have for some positive constants $c_i$
\[\left|
\int_{{\bf R}^N \setminus B_{1}(\bar{0})} g(\bar{x}) 
e^{- \frac{1}{2}\bar{x}{\cal H}\bar{x}' }d\bar{x}\right|	\leq c_3\int_{{\bf R}^N \setminus B_{1}(\bar{0})} \|\bar{x}\|^k
e^{-\frac{c_2}{a^3}\|\bar{x}\|^2-\frac{1}{4}\bar{x}{\cal H}\bar{x}'
}d\bar{x}\leq c_4^N e^{- \frac{c_2}{2a^3}}\frac{ 
	(2\pi)^{ \frac{N}{2} } 
}{
	\sqrt{\det {\cal H}} 
}\]
\[= \frac{ 
	(2\pi)^{ \frac{N}{2} } 
}{
	\sqrt{\det {\cal H}} 
} \ o\left(e^{- \frac{c_2}{3}
	\lambda^{3/2}
} \right).
\]
Making use of the last bound and (\ref{Ja22}) we 
derive from (\ref{Ja25}) 

\begin{equation}\label{Ja24}
\int_{[0,1]^N }  g(\bar{x}-\bar{a}^{\circ})
e^{- H_{\lambda}^{\circ, N}(\bar{x})  }d\bar{x}
\end{equation}
\[
=e^{-H_{\lambda}^{\circ, N}(\bar{a}^{\circ})
}
\int_{{\bf R}^N
}  
g(\bar{x}) e^{- \frac{1}{2}\bar{x}{\cal H}\bar{x}'    }d\bar{x}\left( 1
+
o\left( \lambda^{-\frac{1}{4}+3q}
\right) \right) \]
\[+ e^{-H_{\lambda}^{\circ, N}(\bar{a}^{\circ})
}\frac{(2\pi)^{\frac{N}{2}}}{\sqrt{\det{\cal H}(N)}}
o\left( e^{-\frac{B}{4}
	({\log \lambda})^2}
\right)\left( 1+O\left(\varepsilon\right) \right)^{\frac{N}{2}},
\]
where $\varepsilon $ is defined by (\ref{Ja20}).

Let us return now to formula (\ref{S13}). Note that relation (\ref{Ja24}) is derived only under assumption (\ref{Ja20l}), i.e., that $\lambda >N$. Hence, replacing $N$ in (\ref{Ja24}) by $2\widetilde{N}$ 
gives us as well
\begin{equation}\label{Ja24n}
\int_{[0,1]^{2\widetilde{N}} }  g(x_1-{a}, \ldots, x_{\widetilde{N}}-a)
e^{- H_{\lambda}^{\circ, 2\widetilde{N}}(x_1, \ldots, x_{2\widetilde{N}})  }d
x_1, \ldots, x_{2\widetilde{N}}
\end{equation}
\[
=e^{-H_{\lambda}^{\circ, 2\widetilde{N}}(\bar{a}^{\circ})
}
\int_{\mathbb{R}^{2\widetilde{N}}}
g(\bar{x}) e^{- \frac{1}{2}\bar{x}{\cal H}(2\widetilde{N})\bar{x}'    }d\bar{x}\left( 1
+
o\left( \lambda^{-\frac{1}{4}+3q}
\right) \right) \]
\[+ e^{-H_{\lambda}^{\circ, 2\widetilde{N}}(\bar{a}^{\circ})
}\frac{(2\pi)^{\frac{2\widetilde{N}}{2}}}{\sqrt{\det{\cal H}(2\widetilde{N})}}
o\left( e^{-\frac{B}{4}
	({\log \lambda})^2}
\right)\left( 1+O\left(\frac{\log \lambda}{\lambda^{\frac{1}{4}}}\right) \right)^{\frac{2\widetilde{N}}{2}},
\]
where vector $\bar{a}^{\circ}$ is $2\widetilde{N}$-dimensional.

Here we shall point out that the desired Gaussian approximation will follow from the last bound only for small enough  $\widetilde{N}$. This explains our choice 
\begin{equation}\label{S3}
\widetilde{N} =\left[N^{\frac{1}{7}} \right]
\end{equation}
as we assume in the theorem.
Then taking  also into account another assumption of the theorem that  $\lambda \sim N^2$, and setting $q= \frac{1}{36}$, we derive from
(\ref{Ja24n})
\begin{equation}\label{S14}
\int_{[0,1]^{2\widetilde{N}} }  g(x_1-{a}, \ldots, x_{\widetilde{N}}-a)
e^{- H_{\lambda}^{\circ, 2\widetilde{N}}(x_1, \ldots, x_{2\widetilde{N}})  }d
x_1, \ldots, x_{2\widetilde{N}}
\end{equation}
\[
= e^{-H_{\lambda}^{\circ, 2\widetilde{N}}(\bar{a}^{\circ})
}\left( 1
+
o\left( N^{-\frac{1}{3}}
\right) \right) \left(
\int_{\mathbb{R}^{2\widetilde{N}}}
g(\bar{x}) e^{- \frac{1}{2}\bar{x}{\cal H}(2\widetilde{N})\bar{x}'    }d\bar{x}+ 
\frac{(2\pi)^{\frac{2\widetilde{N}}{2}}}{\sqrt{\det{\cal H}(2\widetilde{N})}}
o\left( e^{-\frac{B}{4}
	({\log \lambda})^2}
\right) \right).
\]
Substituting  now (\ref{S14}) into  (\ref{S13}) we derive
\begin{equation}\label{S15}
{\bf E} g\left(
{X}^{\circ, N}_{1, \lambda}-a, \ldots, {X}^{\circ, N}_{\widetilde{N}, \lambda}-a\right)
\end{equation}
\[
=\left(1	+
\ O\left(e^{-\alpha \widetilde{N} }
\right)\right) \frac{\int_{[0,1]^{2\widetilde{N}}} g(x_{1}-a,  \ldots ,x_{\widetilde{N}}-a)e^{-H_{\lambda}^{\circ, 2\widetilde{N}}(
		x_{1},  \ldots ,x_{2\widetilde{N}}
		)}
	dx_{1} \ldots d x_{2\widetilde{N}}}{
	\int_{[0,1]^{2\widetilde{N}}} 
	e^{-H_{\lambda}^{\circ, 2\widetilde{N}}(
		x_{1},  \ldots ,x_{2\widetilde{N}}
		)}
	dx_{1} \ldots d x_{2\widetilde{N}}
}
\]
\[
= 
\frac{ \sqrt{\det {\cal H}(2\widetilde{N})}} {
	(2\pi)^{ \frac{2\widetilde{N}}{2} } 
}
\int_{{\bf R}^{2\widetilde{N}}
}  
g(\bar{x}) e^{- \frac{1}{2}\bar{x}{\cal H}(2\widetilde{N})\bar{x}'    }d\bar{x}
\left( 1
+
o\left( N^{-\frac{1}{3}}
\right) \right)
+
o\left( e^{-\frac{B}{4}
	({\log \lambda})^2}
\right),
\]
and hence statement 
(\ref{Ma39}) follows.

\hfill$\Box$

\bigskip

{\bf Proof of Theorem \ref{Lexp1}.} First we compute 
the Hessian ${\cal H}_1(N)$ of $H_{\lambda}$ at ${\bar a}$:
\[
{\small \left(
	\begin{array}{cccccc}
	\frac{2\beta}{a_1^3}
	+\frac{\gamma}{2(a_1+a_2)^3}
	& \frac{\gamma}{4(a_1+a_2)³} & 0 & \ldots&  0 & 0 \\
	\frac{\gamma}{4(a_1+a_2)^3}
	& \frac{2\beta}{a_2³}+\frac{\gamma}{2(a_1+a_2)³}+\frac{\gamma}{2(a_2+a_3)³}& \frac{\gamma}{4(a_2+a_3)³} & 0&  \ldots&  0 \\
	\ldots & & & & & \\
	0 & 0 &  \ldots &   0 & \frac{\gamma}{4(a_{N-1}+a_N)³} &
	\frac{2\beta}{a_N^3}
	+\frac{\gamma}{2(a_{N-1}+a_N)^3}                            
	\end{array}
	\right).}
\]
Observe that ${\cal H}_1$ is non-negative definite due to the fact that $H_{\lambda}$ is convex with the minimum at ${\bar a}$. Furthermore, it follows directly by the 3-diagonal form of ${\cal H}_1$ that its rank is at least $N-1$. Then it is straightforward to check that the first row is not a linear combination of the rest, and therefore ${\cal H}_1$ is positive definite.
Assuming that $\lambda$ is of order $N²$, the entries of vector ${\bar a}$, which is the solution to system (\ref{H2}),  all must be of order $\lambda^{-1/2}$, as for the vector ${\bar a}$. 
Hence, the entries of ${\cal H}_1$ have the same scaling $\lambda^{-3/2}$ with respect to $\lambda$ as in ${\cal H}$. This allows us to use the same approximation technique as in the proof of Theorem \ref{Lexp0}, which leads to the statements of the theorem. \hfill$\Box$
\bigskip

\subsection{Expectations and covariances}\label{EC}

The following proposition (proved below in Appendix) describes 
the inverse of  matrix ${\cal H}$ which is the covariance matrix $\Lambda = {\cal H}^{-1}$ for the Gaussian approximating vector in 
Theorem \ref{Lexp0}.

\begin{prop}\label{Pcov}
	For any $\beta >0$ and $0\leq \gamma \leq \beta$  let
	\[
	\delta =  \frac{\gamma}{4\beta +\gamma+ 2\sqrt{ 4\beta^2+2\beta \gamma}}
	\]
	(as defined in (\ref{del})). Then $\Lambda(n) = {\cal H}^{-1}(n)$ satisfies the following:
	\begin{equation}\label{L11*}
	\Lambda_{11}(n)= \left( \frac{2\beta + \gamma}{2\lambda}\right)^{3/2} \frac{1}{2\beta + \gamma (1 - \delta)/2}+O(e^{-\alpha n}),
	\end{equation}
	and for all $1\leq k\leq [\frac{n+1}{2}]$
	\begin{equation}\label{23J71*}
	\Lambda_{1k}(n) =
	\Lambda_{11}(n)
	(-\delta)^{k-1}+O(e^{-\alpha n})
	\end{equation}
	\[=
	\left( \frac{2\beta + \gamma}{2\lambda}\right)^{3/2} \frac{1}{2\beta + \gamma (1 - \delta)/2}
	(-\delta)^{k-1}+O(e^{-\alpha n}),\]
	where $\alpha$ is some positive constant. \hfill$\Box$
\end{prop}

\begin{rem}\label{remM}
	Observe that by the symmetry it holds that for all $1\leq k \leq N$
	\begin{equation}\label{Msyn}
	{\bf E}X_{k,\lambda}^{\circ}
	={\bf E}X_{1,\lambda}^{\circ},
	\end{equation}
	and for all $i\leq j \leq N$
	\begin{equation}\label{Covsym}
	{\bf Cov} (X_{i,\lambda}^{\circ},  X_{j,\lambda}^{\circ})={\bf Cov} (X_{1,\lambda}^{\circ},  X_{1+j-i,\, \lambda}^{\circ}).
	\end{equation}
\end{rem}

As an immediate corollary of this remark and Theorem \ref{Lexp0} we get the following result.

\begin{cor}\label{ExCov}
	For all $1\leq i \leq N$ and $\lambda \sim N^2$ we have
	\begin{equation}\label{Ja8*}
	{\bf E}X_{i,\lambda}^{\circ}
	= a\left(1+ o\left(\lambda^{-\frac{1}{6}}
	\right)\right),
	\end{equation}
	and if $|i-j|\leq N^{1/7}$
	\begin{equation}\label{M}
	{\bf Cov} (X_{i,\lambda}^{\circ},  X_{j,\lambda}^{\circ})= ({\cal H}^{-1})_{ij}\left(1+ o\left(\lambda^{-\frac{1}{6}}
	\right)\right)+o\left(e^{-(\log \lambda)^{\frac{3}{2}}}
	\right),
	\end{equation}
	where $({\cal H}^{-1})_{ij}$ are the elements of the matrix inverse to ${\cal H}(n)$ with $n=2\left[N^{1/7}\right]$.  
	All the remaining terms are uniform in $i,j=1, \ldots, N.$
\end{cor}
\hfill$\Box$

Similarly Theorem \ref{Lexp1} yields the following result.
\begin{cor}\label{ExCov1}
	For all $1\leq i\leq N$ and $\lambda \sim N^2$ we have
	\begin{equation}\label{Ja8*1}
	{\bf E}X_{i,\lambda}
	= a_i\left(1+ o\left(\lambda^{-\frac{1}{6}}
	\right)\right),
	\end{equation}
	and if $|i-j|\leq N^{1/7}$
	\begin{equation}\label{M1}
	{\bf Cov} (X_{i,\lambda},  X_{j,\lambda})= ({\cal H}_1^{-1})_{ij}\left(1+ o\left(\lambda^{-\frac{1}{6}}
	\right)\right)+o\left(e^{-(\log \lambda)^{\frac{3}{2}}}
	\right),
	\end{equation}
	where $({\cal H}_1^{-1})_{ij}$ are the elements of the matrix inverse to ${\cal H}_1$.  
	All the remaining terms are uniform in $i,j=1, \ldots, N.$
\end{cor}
\hfill$\Box$

Proposition \ref{Pcov} helps us  to compute the covariances  in Corollary \ref{ExCov}.

\begin{cor}\label{CorMay1} 
	Let $\lambda \sim N^2$.
	Then 
	\begin{equation}\label{May19}
	{\bf Var} (X_{k,\lambda}^{\circ})
	= \left( \frac{2\beta + \gamma}{2\lambda}\right)^{3/2} \frac{1}{2\beta + \gamma (1 - \delta)/2}\left(1+ o\left(\lambda^{-\frac{1}{6}}
	\right)\right)+o\left(e^{-(\log \lambda)^{\frac{3}{2}}}
	\right),
	\end{equation}	
	and for all $2\leq k \leq \left[
	\frac{N+1}{2} \right]$
	\begin{equation}\label{May21}
	{\bf Cov} (X_{1,\lambda}^{\circ},  X_{k,\lambda}^{\circ})
	= 	{\bf Cov} (X_{1,\lambda}^{\circ},  X_{N-k+2,\lambda}^{\circ})
	\end{equation}
	\[        =(-1)^{k-1}
	\left( \frac{2\beta + \gamma}{2\lambda}\right)^{3/2} 
	\frac{\delta^{k-1}}{2\beta + \gamma (1 - \delta)/2}
	\left(1+ o\left(\lambda^{-\frac{1}{6}}
	\right)\right)+o\left(e^{-(\log \lambda)^{\frac{3}{2}}}
	\right)
	.\]
\end{cor}

\noindent{\bf Proof.}
By the symmetry we have in (\ref{M})
\[ ({\cal H}^{-1})_{ij}= \Lambda _{1 \, (|i-j|+1)}(n)\]
where $n=2\left[ N^{1/7}\right]$. Therefore substituting  formula (\ref{23J71*}) (with $n=2\left[ N^{1/7}\right]$) into  (\ref{M}) we
get for all $0\leq k \leq N^{1/7}$
\begin{equation}\label{M1*}
{\bf Cov} (X_{1,\lambda}^{\circ},  X_{1+k,\lambda}^{\circ})= 
\Lambda _{1 \, (k+1)}(n)
\left(1+ o\left(\lambda^{-\frac{1}{6}}
\right)\right)+o\left(e^{-(\log \lambda)^{\frac{3}{2}}}
\right)
\end{equation}
\[=\left( 
\left( \frac{2\beta + \gamma}{2\lambda}\right)^{3/2} \frac{1}{2\beta + \gamma (1 - \delta)/2}
(-\delta)^{k }+O\left(e^{-\alpha N^{1/7}/2}\right)\right)\left(1+ o\left(\lambda^{-\frac{1}{6}}
\right)\right)\]
\[+o\left(e^{-(\log \lambda)^{\frac{3}{2}}}
\right)
\]
\[= 
\left( \frac{2\beta + \gamma}{2\lambda}\right)^{3/2} \frac{(-\delta)^{k }}{2\beta + \gamma (1 - \delta)/2}
\left(1+ o\left(\lambda^{-\frac{1}{6}}
\right)\right)+o\left(e^{-(\log \lambda)^{\frac{3}{2}}}
\right).
\]
This together with assumption that $\lambda \sim N^2$ implies 
(\ref{May19}) when $k=0$  in (\ref{M1*}), and (\ref{May21}) 
when $1\leq k \leq N^{1/7}$ in (\ref{M1*}).

In the remaining case when $N^{1/7}<k<\frac{N}{2}$
we simply use result (\ref{Ea21}) on the exponential decay of correlations
\begin{equation}\label{M1n}
{\bf Cov}  (X_{1,\lambda}^{\circ},  X_{1+k,\lambda}^{\circ})=O\left(e^{-\alpha ' k}\right)
\end{equation}
for some positive $\alpha '$ uniformly in $\lambda$.
Under assumption that   $\lambda \sim N^2$ the  asymptotic (\ref{May21}) is still valid
as well for $N^{1/7}<k<\frac{N}{2}$,  since in this case by (\ref{M1n}) 
$${\bf Cov}  (X_{1,\lambda}^{\circ},  X_{1+k,\lambda}^{\circ})= O\left(e^{-\alpha ' k}\right)  = 
o\left(e^{-(\log \lambda)^{\frac{3}{2}}}
\right),$$
as in (\ref{May21}).
\hfill$\Box$

\medskip

Corollary \ref{CorMay1} reveals an interesting covariance structure of the distribution of location of particles: two interspacings are positively correlated if the number of spacings between them is odd, otherwise, they are negatively correlated. Nevertheless, as 
it follows directly by the results of Corollary \ref{CorMay1},
the covariance between any spacing and the total sum of the spacings remains to be positive for all parameters $\beta>0$ and $0\leq  \gamma\leq \beta.$ This yields as well a strictly positive variance of $\sum_{k=1}^NX_{k,\lambda}^{\circ}$ as described below.

\begin{cor}\label{CorMay2} 
	For all positive parameters 
	we have 
	\begin{equation}\label{23Ja8}
	{\bf Var} \left(\sum_{k=1}^NX_{k,\lambda}^{\circ}\right)=
	N 	{\bf Cov} \left(X_{1,\lambda}^{\circ},\sum_{k=1}^NX_{k,\lambda}^{\circ}\right)
	\end{equation}
	\[= \left(\frac{2\beta + \gamma}{2\lambda} \right)^{\frac{3}{2}} N \frac{1}{2\beta + \gamma (1 - \delta)/2} \ \frac{1-\delta}{1+\delta}\left(1+ o\left(\lambda^{-\frac{1}{6}}
	\right)\right),\]
	where the principal term is strictly positive.
\end{cor}
\hfill$\Box$

\subsection{Setting value for the Lagrange multiplier}\label{VLM}

Now we can choose  $\lambda=\lambda(N)$ which solves (\ref{ES}). By Corollary \ref{ExCov1} we have
\[{\bf E}X_{k,\lambda}
= a_k\left(1+ o\left(\lambda^{-\frac{1}{5}}
\right)\right),\]
while by (\ref{De22n5}) and Corollary \ref{ExCov} we get
for all $1\leq k \leq N$ and $\lambda \sim N^2$ 
\begin{equation}\label{Nov6}
{\bf E}X_{k,\lambda}
= \left(1+o\left(e^{-\alpha k}\right)\right){\bf E}X_{k,\lambda}^{\circ}=
\left(1+o\left(e^{-\alpha k}\right)\right)\left(1+ o\left(\lambda^{-\frac{1}{5}}
\right)\right)a.
\end{equation}	
Combination of these results gives us
\[
\frac{a_k}{a}-1 = o\left(N^{-\frac{2}{5}}
\right)+o\left(e^{-\alpha k}\right),\]
which confirms conditions of Lemma \ref{PA1} when $N$ is large. Therefore Lemma \ref{PA1} implies 
\begin{equation}\label{ak}
{a_k}={a} +O \left({\eta^k}\right),
\end{equation}
where
\[ \eta = 1+4\frac{\beta}{\gamma} -\sqrt{16\left(\frac{\beta}{\gamma} \right)^2 +8\frac{\beta}{\gamma}}<1.\]

Now we can rewrite 
(\ref{ES}) as
\begin{equation}\label{ES1}
1=\sum_{k=1}^N \mathbb{E} X_{k,\lambda} =
\sum_{k=1}^{\sqrt N} 
O\left(\frac{1}{N} \right)
+
\sum_{k={\sqrt N} } ^N
\left( {a} +O \left({\eta^k}\right)\right)
=\sqrt{\frac{2\beta + \gamma}{2\lambda}}N \left(1+ O(N^{-1/2})\right),
\end{equation}
as $N\rightarrow \infty$.
Hence, the solution to (\ref{ES}) is
\begin{equation}\label{lambda}
\lambda=\lambda(N)=  \sqrt{\frac{2}{2\beta + \gamma}}\ N^2
\left(1+ O(N^{-1/2})\right),
\end{equation}
which implies 
\begin{equation}\label{fa}
a=\frac{1}{N}\left(1+ O(N^{-1/2})\right).
\end{equation}
Note also that
with this choice of $\lambda$ we get from (\ref{23Ja8}) setting  $x:= \gamma/\beta$
\begin{equation}\label{May31}
{\bf Var} \left(\sum_{k=1}^N X_{k,\lambda}^{\circ}\right)   =  \frac{1}{\beta N^{2} }\,  g(x) \left(1+ o\left(N^{-\frac{1}{6}}
\right)\right),
\end{equation}
where
\[ g(x)
=\frac{2}{
	4 + \frac{ 4x +2x\sqrt{ 4+2x}}{4 +x+ 2\sqrt{ 4+2x}}
} \ \frac{2 + \sqrt{ 4+2x}}{2 +x+ \sqrt{ 4+2x}}>0.  \]

\begin{rem}\label{Rg}
	Function $g(x)$	is monotone decreasing on $R^+$ with the maximal value $g(0)=\frac{1}{2}$, i.e., when $\gamma=0$.
\end{rem}

To prove the statements of Theorem \ref{T1}
we shall consider  conditional expectation of $X^{\circ}_{k, \lambda}$ given 
$\sum_{k=1}^N X_{k, \lambda}^{\circ}=1.$ 
Observe that (\ref{Nov6}) and (\ref{ES1}) imply as well that for some value
\[
\lambda= \sqrt{\frac{2}{2\beta + \gamma}}\ N^2
\left(1+ O(N^{-1/2})\right)
\]
(which differs from the solution (\ref{lambda}) to (\ref{ES}) only in the small term $O(N^{-1/2})$), we have
\begin{equation}\label{Nov7}
\sum_{k=1}^N \mathbb{E} X_{k, \lambda}^{\circ}=1.
\end{equation}
From now on we assume (\ref{Nov7}) to be held, implying
by the symmetry
\begin{equation}\label{Nov8}
m^{\circ}:=\mathbb{E} X_{k, \lambda}^{\circ}=\frac{1}{N}.
\end{equation}
Furthermore, by the symmetry again it follows (moreover for any value of $\lambda$) that 
\begin{equation}\label{T1E}
\mathbb{E}\left\{ {X}_{k,\lambda}^{\circ}
\mid 
\sum_{k=1}^N X_{k, \lambda}^{\circ}=1\right\}	= 
\frac{1}{N}.
\end{equation}

Next we shall study conditional moments of the second order.

\subsection{Central Limit Theorem}\label{CLT}
Let $f_{S}$ and 
$f_{X_k^{\circ},S}$ denote
the densities of $S$ and of $(X_k^{\circ}, S)$, respectively, and let
$\phi_{S}$ and 
$\phi_{X_k^{\circ},S}$ denote
the corresponding characteristic functions.
We study conditional expectation 
\begin{equation}\label{23Ja4}
\mathbb{E}\left\{ \left( {X}_{k,\lambda}^{\circ}\right)^2
\mid 
\sum_{k=1}^N X_{k, \lambda}^{\circ}=1\right\}	
= \frac{
	\int_{0}^1 x^2f_{X_1^{\circ},S}(x,1-x)
	d{x}
}{
	f_{S^{\circ}_{N,\lambda}}(1)
},
\end{equation}
where the equality follows by the symmetry discussed above.
This instructs us to study the densities involved in the last formula. Therefore we begin with establishing an analogue of Central Limit Theorem. 

Let $N\times N$ covariance matrix $\Lambda={\cal H}^{-1}(N)$ be as defined in Theorem \ref{Lexp0} and described in Proposition \ref{Pcov} .

\begin{lem}\label{TJ1}(Central Limit Theorem)
	Assume that (\ref{Nov7}) holds.
	Let $\bar{ \cal Z} \sim{\cal N}(\bar{0}, \Lambda)$.
	Set
	$$\sigma_N^2:= \mbox{{\bf Var}}\left(	\sum_{k=1}^N {\cal Z}_k\right),$$
	and 
	$$\tilde{X_k}=\frac{X^{\circ}_{k, \lambda}-\frac{1}{N}}{\sigma_N}, \ \ \ \ 
	\tilde{S}=\frac{\sum_{k=1}^N \left(X^{\circ}_{k, \lambda}-\frac{1}{N}\right)}{\sigma_N}, \ \ \ \tilde{{\cal Z}}_k=\frac{{\cal Z}_k}{\sigma_N},\ \ \ \Sigma = \frac{\sum_{k=1}^N{\cal Z}_k}{\sigma_N}.$$ 
	Then under the assumptions of Theorem \ref{Lexp0} 
	we have for all $x\in \mathbb{R}$
	\begin{equation}\label{23C1}
	f _{\tilde{S}}(x)=f_{\Sigma}(x)+ O\left(N^{-1/6}  \right),
	\end{equation}
	where the error term in uniform in $x$.
\end{lem}

\noindent
{\bf Proof.} Observe that by the assumption $\lambda$ is a function of $N$, hence $N$ will be the only variable below while $\beta>0$
and $\gamma\geq 0$ are fixed arbitrarily.

The proof is largely based on the ideas and results of
\cite{S}. In our setting the random variables $X^{\circ}_{k, \lambda}$
are indexed by $k\geq 1$, which corresponds to the 
case of random variables in \cite{S} on a one-dimensional lattice
(i.e., $d=1$ in notation of \cite{S}).

The crucial conditions for the proof of CLT as in \cite{S} 
are guaranteed in our case by (\ref{Nov3}) and by Corollary \ref{CorMay1}. The latter implies similar to the proof of  (\ref{23Ja8}), that there is a constant $A_1 $ (independent of $N$) such for any set $I \subseteq \{1, \ldots, N\}$ of consecutive indices 
\begin{equation}\label{Nov1}
Var \left(\sum_{i\in I}X^{\circ}_{k, \lambda}\right)\geq \frac{A_1}{N^3}|I|.
\end{equation}
Properties (\ref{Nov3}) and (\ref{Nov1}) allow us to
follow closely the lines of \cite{S}  and derive for any set $I \subseteq \{1, \ldots, N\} $ of consecutive indices
\begin{equation}\label{Nov5}
\mathbb{E} \left(
\sum_{k\in I} \left(X^{\circ}_{k, \lambda}
-\mathbb{E}X^{\circ}_{k, \lambda}\right)
\right)^4
\leq A_2 \left( \frac{1}{N^3} \right)^2|I|^2.
\end{equation}
for some constant $A_2$.
Observe, however, that our bounds on the right both in (\ref{Nov1}) and (\ref{Nov5}) depend on $N$ unlike the case in \cite{S}. This requires some extra care as we clarify below.

For any random vector $\xi \in \mathbb{R}^k$ let us denote its characteristic function by 
\[\phi_{\xi}(\bar{t})=\mathbb{E}e^{i(\bar{t}, \xi)}.\]
Making use of 
Fourier inverse formula for the densities 
we derive first for an arbitrarily fixed $T>0$
\begin{equation}\label{I1}
\left|f_{\tilde{S}}(x)-f_{\Sigma}(x)\right|
=  \left|\frac{1}{{2\pi} }
\int_{-\infty}^{\infty}e^{-itx}\left( \phi_{\Sigma}(t)-\phi_{\tilde{S}}(t)\right)
dt\right|
\end{equation}
\[
\leq\frac{1}{{2\pi} }
\int_{|t|\leq T}
\left| \phi_{\Sigma}(t)-\phi_{\tilde{S}}(t)\right|
dt + \frac{1}{{2\pi} } \left|\int_{|t|>T} e^{-itx}
\phi_{\tilde{S}}(t)dt\right|+o\left(e^{-T^2/3}\right),
\]
where $\phi_{\Sigma}(t)=e^{-t^2/2}$.
Let us now explain briefly
the idea of \cite{S} 
how  to obtain an upper bound on the difference in the first integral on the right. The sum $\tilde{S}$ is split into partial sums over the sets of indices well separated 
($W_i$ below) and the sums over separating indices  ($R_i$).
Then with a proper choices of set sizes the sum $\tilde{S}$ can be approximated by a sum of almost independent terms. To execute this idea consider a decomposition
of $\{1, \ldots, N \}$ into subsets of consecutive indices
\[
\{1, \ldots, N \}= 
\{1, \ldots, p\} 
\cup \{p+1, \ldots, p+q\} \cup \ldots \]
\[=:W_1\cup R_1 \cup W_2 \cup R_2 \cup \ldots \cup W_k \cup R_k,\]
where it is assumed that 
$q=q(N) \rightarrow \infty$,  while for some positive $\varepsilon$
\begin{equation}\label{vare}
p=N^{\varepsilon}q =o(N), \ \ \ k=\left[ \frac{N}{p+q}\right] \rightarrow \infty, \ \ \ \mbox{as} \ \ \  N \rightarrow \infty, 
\end{equation}
so that
\[|W_i|=p, i=1, \ldots, k, \ \ \ |R_j|=q, j=1, \ldots , k-1, \ \ \ |R_k|=q+N-(p+q)\left[ \frac{N}{p+q}\right].\]
Then a number of the following approximation (in distribution) steps are taken as explained in \cite{S}:
\[\tilde{S}= \frac{\sum_{k=1}^NX^{\circ}_{k, \lambda}-1}{\sigma_N} \approx
\frac{\sum_{i=1}^k \sum_{j\in W_i} \left( X^{\circ}_{j, \lambda}-\frac{1}{N}\right)}{\sigma_N} \approx
\frac{\sum_{i=1}^k \zeta_i}{\sqrt{k Var(\zeta_1)}} \approx
\frac{\sum_{i=1}^k \zeta_i'}{\sqrt{k Var(\zeta_1)}} \approx \Sigma,
\]
where $$\zeta _i= \sum_{j\in W_i} \left( X^{\circ}_{j, \lambda}-\frac{1}{N}\right)$$
(recall that $\mathbb{E}X^{\circ}_{j,\lambda}=\frac{1}{N}$),
while
$\zeta _i', i\geq 1, $ are {\it independent} random variables  with the same distribution:  
\[\zeta _i'\stackrel{d}{=}\zeta _i.\]
Note that by this definition all
$\zeta _i$,
$\zeta _i'$, $i\geq 1$, are identically distributed.

With the arguments from \cite{S} the above chain of approximations leads in our case to 
\begin{equation}\label{Nov10}
\left| \phi_{\tilde{S}}(t) - \phi_{\Sigma}(t)\right|= \left| {\mathbb E}e^{it\tilde{S}} - e^{-t^2/2}\right|\leq O\left(N^{-\varepsilon/2}\right)t+O\left(e^{-\alpha q/2}\right) +O\left(t^4e^{-t^2/4}\right)\frac{1}{k},
\end{equation}
as long as  
\begin{equation}\label{Nov9}
|t|\leq T=A_3 \sqrt{k},
\end{equation}
for   an arbitrarily fixed positive constant $A_3$. Here $\alpha>0$ is the same as in  (\ref{Nov3}), and $\varepsilon$ and $k$ are defined in (\ref{vare}).
Hence, when
\begin{equation}\label{Novq}
q\geq \log^2 N,
\end{equation}
bound (\ref{Nov10}) yields
\begin{equation}\label{Nov11}
\int_{|t|\leq T}
\left| \phi_{\Sigma}(t)-\phi_{\tilde{S}}(t)\right|
dt\leq O\left(N^{-\varepsilon/2}\right)T^2+O\left(\frac{1}{k}\right)
=O\left(N^{-\varepsilon/2}k+\frac{1}{k}\right),
\end{equation}
where the last equality is due to the definition of $T$ in (\ref{Nov9}).
Let us choose now an optimal $\varepsilon$ to minimize the right hand side in
(\ref{Nov11}).
Recall that 
by definition (\ref{vare})
\begin{equation}\label{Novp}
O\left(N^{-\varepsilon/2}k+\frac{1}{k}\right)= O\left(N^{-\varepsilon/2} \frac{N}{p}+\frac{p}{N}\right),
\end{equation}
which is of minimal order in $N$ when 
\begin{equation}\label{Novp1}
p= N^{1-\frac{\varepsilon}{4}}. 
\end{equation}
Hence,  if $p$ is chosen as in (\ref{Novp1}) with  $\varepsilon$ such that  
\begin{equation}\label{Novpq}
0<{\varepsilon}<1-\frac{\varepsilon}{4},
\end{equation}
both conditions (\ref{vare}) and (\ref{Novq}) are satisfied. 
Setting now
\begin{equation}\label{Nov12}
T=A_4 N^{\frac{\varepsilon}{8}} \leq A_3 \sqrt{k}
\end{equation}
for some positive constant $A_4,$
we derive from (\ref{Nov11}) and (\ref{Novp}) with a help of (\ref{Novpq}) that
\begin{equation}\label{Nov12*}
\int_{|t|\leq T}
\left| \phi_{\Sigma}(t)-\phi_{\tilde{S}}(t)\right|
dt\leq O\left(N^{-\frac{\varepsilon}{4}}\right)
\end{equation}
for any 
\[0<\varepsilon <\frac{4}{5}.\]
From now on we set 
\begin{equation}\label{NovEps}
\varepsilon = \frac{2}{3}, \ \ \ T=A_4 N^{\frac{1}{12}}.
\end{equation}

The upper bound for the remaining integral over $|t|>T$ in (\ref{I1})
is complicated (unlike in \cite{S}) by the dependence on $N$ of the parameters of distribution of 
our variables $\bar{X}^{\circ}_{\lambda}$.

We shall first find an upper bound for $\phi_{\tilde{S}}(t)$ for $|t|>T$ making
use of the results on the local Gaussian approximation established in Theorem \ref{Lexp0}. 

Taking into account that in ${H}_{\lambda}^{\circ}$
the variables $x_k$ and $x_n$ are, loosely speaking, "independent" if $|n-k|>2$,
let us split ${\tilde{S}}$ again as follows
\begin{equation}\label{Nov15}
{\tilde{S}} =\sum_{k=1}^N {\tilde{X}_k} 
=: \sum_{n=0}^{j-1} \sum_{k=2}^{{\tilde N}-2}  {\tilde{X}}_{n{\tilde N}+k}
+  \left( {\tilde{X}_{N}} + {\tilde{X}}_{1}+
\sum_{n=0}^{j-1} \sum_{k=-1}^{1} {\tilde{X}}_{n{\tilde N}+k}
\right) + {\cal R}
\end{equation}
where ${\tilde N}=\left[N^{\frac{1}{7}}\right]$ as defined  in Theorem \ref{Lexp0},
and 
\begin{equation}\label{j}
3<j\ll N^{1-\frac{1}{7}}
\end{equation}
is to be fixed  later.
The first sum on the right in (\ref{Nov15}) consists of partial sums of consecutive variables separated by 3 indices, next bracketed term contains all the boundary variables for the  sums in the first term, and the last term ${\cal R}$ contains all the rest.
Denote 
\begin{equation}\label{Nov16}
\xi_n =\sum_{k=2}^{{\tilde N}-2}  {\tilde{X}}_{n{\tilde N}+k},
\end{equation}
and let
\[B:=\{1,2\} \cup \{n{\tilde N}+k: k=-1, 0, 1, n=1, \ldots, j-1\}\]
denote the set of indices for the variables in the brackets in (\ref{Nov15}).
Then by (\ref{Nov15})
we have
\begin{equation}\label{Nov17}
{\tilde{S}} =\sum_{n=0}^{j-1}
\xi_n
+  
\sum_{k\in B} {\tilde{X}}_{k}
+ {\cal R}.
\end{equation}
Observe that due to the form of distribution (\ref{H0c}) conditionally on the values ${\tilde{X}}_{k}, k\in B$ (or, equivalently, on $X^{\circ}_{k, \lambda}, k\in B$)
the random variables ${\cal R}, \xi_0, \ldots , \xi_{j-1}$, are independent. Note also that
$\xi_0, \ldots , \xi_{j-1}$ are 
identically distributed.
These properties allow us to derive
\begin{equation}\label{Nov18}
\left| {\mathbb E}e^{it\tilde{S}} \right| \leq \left| {\mathbb E}e^{it \left(\sum_{n=0}^{j-1}
	\xi_n
	+  
	\sum_{k\in B} {\tilde{X}}_{k}
	+ {\cal R}\right)}{\bf I}\left\{ \max _{k\in B} \left|X^{\circ}_{k, \lambda}-\frac{1}{N}\right|\leq 
\left(\frac{1}{N} \right)^{5/4}
\right\}\right| 
\end{equation}
\[+{\bf P}\left\{ \max _{k\in B} \left|X^{\circ}_{k, \lambda}-\frac{1}{N}\right|>
\left(\frac{1}{N} \right)^{5/4}
\right\} \]
\[
\leq
\left(  
\max_{\left|x_k\right|\leq 
	{N}^{-5/4}, \ k\in B
}\left|
{\mathbb E}
\left\{
e^{it
	\left(
	\tilde{X}_{N}+
	\tilde{X}_{1}+\xi_0
	+ 
	\tilde{X}_{ \tilde {N}-1}+
	\tilde{X}_{\tilde {N}}
	\right)} \mid   X^{\circ}_{k, \lambda}-\frac{1}{N}=x_k, k\in B
\right\}\right|
\right)^j
\]
\[+ |B| {\bf P}\left\{ \left|X^{\circ}_{1, \lambda}-\frac{1}{N}\right|>
\left(\frac{1}{N} \right)^{5/4}
\right\}.\]
Consider in the last formula
\[\phi_{\xi, \bar{x}}(t):= {\mathbb E}
\left\{
e^{it
	\left(
	\tilde{X}_{N}+
	\tilde{X}_{1}+\xi_0
	+ 
	\tilde{X}_{ \tilde {N}-1}+
	\tilde{X}_{\tilde {N}}
	\right)} \mid   X^{\circ}_{k, \lambda}-\frac{1}{N}=x_k, k\in B
\right\}\]
\[=
{\mathbb E}
\left\{
e^{it \left( \frac{ X^{\circ}_{N, \lambda}-\frac{1}{N}}{\sigma_N} + \sum_{k=1}^{\tilde {N}}
	\frac{ X^{\circ}_{k, \lambda}-\frac{1}{N}}{\sigma_N}\right)
}
\mid   X^{\circ}_{k, \lambda}-\frac{1}{N}=x_k, k = N,1, \tilde {N}-1, \tilde {N}
\right\}.
\]
Let us write ${\cal Z}_{-1}={\cal Z}_{ {N}}$. By the Theorem \ref{Lexp0}
\begin{equation}\label{Nov20} 
\phi_{\xi, \bar{x}}(t)=
{\mathbb E}
\left\{
e^{i\frac{t }{\sigma_N} \sum_{k=-1}^{\tilde {N}}{\cal Z}_k 
}
\mid   {\cal Z}_k=x_k, k = -1,1, \tilde {N}-1, \tilde {N}
\right\}
\left(1	+
\ o\left(N^{-\frac{1}{3}}
\right)\right)
\end{equation}
\[
+o\left( e^{-
	({\log N})^{3/2}}
\right).\]
By the properties of Normal distribution the conditional  distribution of
\[\sum_{k=-1}^{\tilde {N}}{\cal Z}_k 
\mid  \left( {\cal Z}_k=x_k, k = -1,1, \tilde {N}-1, \tilde {N}\right)\]
is again Normal with the parameters which are straightforward to compute: the mean is  $$O\left( |x_{-1}|+|x_1|+|x_{\tilde {N}-1}|+|x_{\tilde {N}}|\right)=O(N^{-5/4}),$$
when each $|x_i|\leq N^{-5/4},$ and variance 
\[(C+o(1))\frac{\tilde {N}}{N}\sigma_N^2 +O\left(N^{-3} \right)\]
for some positive constant $C$. Hence, for all large $N$ 
\begin{equation}\label{Nov21}
|\phi_{\xi, \bar{x}}(t)|\leq e^{-\frac{t^2}{3\sigma_N^2}
	C\frac{\tilde {N}}{N}\sigma_N^2
}
\left(1	+
\ o\left(N^{-\frac{1}{5}}
\right)\right)+o\left( e^{-
	({\log N})^{3/2}}
\right)\leq e^{-\frac{t^2}{4}
	C\frac{\tilde {N}}{N}
}+o\left( e^{-
	({\log N})^{3/2}}
\right)
\end{equation}
uniformly in $\bar{x}$ when all $|x_i|\leq N^{-5/4}.$
Making use of this bound in (\ref{Nov18}) we get 
\begin{equation}\label{Nov22}
\left| {\mathbb E}e^{it\tilde{S}} \right| \leq
\left(  e^{-\frac{t^2}{4}
	C\frac{\tilde {N}}{N}
}+o\left( e^{-
	({\log N})^{3/2}}
\right)
\right)^j
+ 3j {\bf P}\left\{ \left|X^{\circ}_{1, \lambda}-\frac{1}{N}\right|>
\left(\frac{1}{N} \right)^{5/4}
\right\}.
\end{equation}
The Gaussian approximation for $X^{\circ}_{1, \lambda}$
established in Theorem \ref{Lexp0} gives us a bound
\begin{equation}\label{Nov23}
{\bf P}\left\{ \left|X^{\circ}_{1, \lambda}-\frac{1}{N}\right|>
\left(\frac{1}{N} \right)^{5/4}
\right\} <O\left(
e^{-c\left(\frac{1}{N} \right)^{5/2} N^3
}\right)
+o\left( 
e^{-({\log N})^{3/2}}
\right)
\end{equation}
for some positive constant $c$.
Hence, for all 
\[t \geq T=A_4 N^{\frac{1}{12}}
\]
we derive from (\ref{Nov22})
\begin{equation}\label{Nov24}
\left| {\mathbb E}e^{it\tilde{S}} \right| < 
O\left(  e^{-
	jN^{\frac{1}{7}+\frac{1}{7}-1}}
\right)
+o\left(j 
e^{-({\log N})^{3/2}}\right).
\end{equation}
Setting in the last formula
\[j= N^{1-\frac{1}{4}},\]
which satisfies (\ref{j}), we get
\begin{equation}\label{Nov25}
\left| \phi_{\tilde{S}}(t) \right|= \left| {\mathbb E}e^{it\tilde{S}} \right| <
e^{-\frac{1}{2}({\log N})^{3/2}}
\end{equation}
for all 
$t \geq T=A_4 N^{\frac{1}{12}}$.
With a help of the last bound, as well as taking into account (\ref{Nov12*})
together with (\ref{NovEps}) we derive from (\ref{I1})
\begin{equation}\label{J1}
\left|f_{\tilde{S}}(x)-f_{\Sigma}(x)\right|
=  \left|\frac{1}{{2\pi} }
\int_{-\infty}^{\infty}e^{-itx}\left( \phi_{\Sigma}(t)-\phi_{\tilde{S}}(t)\right)
dt\right|
\end{equation}
\[
\leq
O\left(N^{-\frac{\varepsilon}{4}}\right)
+ \frac{1}{{2\pi} } \left|\int_{T}^{N^6} e^{-itx}
\phi_{\tilde{S}}(t)dt\right|
+ \frac{1}{{2\pi} } \left|\int_{|t|>N^6} e^{-itx}
\phi_{\tilde{S}}(t)dt\right|+o\left(e^{-T^2/3}\right),
\]
\[=
O\left(N^{-\frac{1}{6}}\right)
+ \frac{1}{{2\pi} } \left|\int_{|t|>N^6} e^{-itx}
\phi_{\tilde{S}}(t)dt\right|.
\]

Consider  the last integral in (\ref{J1}).
Recall first that
\begin{equation}\label{23J25}
f_{{S}}(x) = \int_{0}^x\int_{0}^{y_{N-1}}\ldots \int_{0}^{y_2}
f_{\bar{X}^{\circ}_{\lambda}}(y_1, y_2-y_1, \ldots, y_{N-1}-y_{N-2}, x-y_{N-1})dy_1 \ldots dy_{N-1},
\end{equation}
with function $ f_{\bar{X}^{\circ}_{\lambda}}$  defined by (\ref{Mt7}) and (\ref{H0c}). Taking into account that
\[f_{\tilde{S}}(x) =\sigma_N f_{{S}}(\sigma_Nx +1), \]
 and
\begin{equation}\label{23J19}
\phi_{\tilde{S}}(t) =
\frac{1}{(it)^2}
\int_{\mathbb{R}}e^{itx}f''_{\tilde{S}}(x)dx,
\end{equation}
we derive
\begin{equation}\label{23J24}
\left| \int_{|t|>N^6} e^{-itx}
\phi_{\tilde{S}}(t) dt \right|\leq  \int_{|t|>N^6} \frac{1}{t^2}dt 
\int_{\mathbb{R}}|f''_{\tilde{S}}(x)|dx,
\end{equation}
where
\begin{equation}\label{23J21}
\int_{\mathbb{R}}\left|f''_{\tilde{S}}(x)\right|dx= {\sigma_N^2}\int_{0}^{{N}}\left|f''_{{S}}(x)\right|dx.
\end{equation}
By a straightforward application of formulas (\ref{Mt7}) and (\ref{H0c})
we get that function $f''_{{S}}(x)$ changes sign at most a finite number (say $L$) of times
on its domain, i.e., $[0,N]$.
Hence,
\begin{equation}\label{23J22}
\int_{0}^{{N}}\left|f''_{{S}}(x)\right|dx \leq 2L\max_{x\in [0,N]} \left|f'_{{S}}(x)\right| .
\end{equation}
Combining   (\ref{23J22}) and   (\ref{23J21}) in  (\ref{23J24})
we obtain 
\begin{equation}\label{23J24n}
\left| \int_{|t|>N^6} e^{-itx}
\phi_{\tilde{S}}(t) dt \right|\leq  \frac{2}{N^{6}}
2L \sigma_N^2
\max_{x} |f_{{S}}'(x) | .
\end{equation}
Here
\begin{equation}\label{23J26}
\max_{x\in [0,N]} \left|f'_{{S}}(x)\right| \leq O(N^2) \max_{x\in [0,N]} f_{{S}}(x),
\end{equation}
which is again a result of straightforward computation based on
(\ref{Mt7}), (\ref{H0c}) and (\ref{23J25}),
and assumption that $\lambda \sim N^2$. Hence, (\ref{23J24n}) yields
\[
\left| \int_{|t|>N^6} e^{-itx}
\phi_{\tilde{S}}(t) dt \right|
\leq N^{-6} 
{\sigma_N^2}O(N^2) \max_{x} f_{{S}}(x)
\]
\begin{equation}\label{23J31}
={O(N^{-5} )}\max_{x} f_{\tilde{S}}(x),
\end{equation}
where we also took into account (\ref{May31}), i.e., that $\sigma_N^2 \sim N^{-2}$.

Equipped with bound (\ref{23J31}) we derive from (\ref{J1})
\begin{equation}\label{23J33}
f_{\tilde{S}}(x) = f_{\Sigma}(x) + 
O\left(N^{-1/6}  \right)
+ O\left(N^{-5}  \right)\max_{x} 
f_{\tilde{S}}(x)
\end{equation}
for all $x$.
Hence, 
$f_{\tilde{S}}(x)$ is uniformly bounded in $x$ as $f_{\Sigma}(x)$ is, and this in turn implies by (\ref{23J33}) that  for all $x $
\begin{equation}\label{23J34}
f_{\tilde{S}}(x) =  f_{\Sigma}(x)+ O\left(N^{-1/6}  \right),
\end{equation}
which is the statement of Lemma. \hfill$\Box$

\begin{cor}\label{CorVa}
	Let $a=1/N$. Then in notation of Lemma \ref{TJ1}
	\begin{equation}\label{CorVar}
	\int_{-\frac{a}{\sigma_N}}^{ \frac{1-a
		}{\sigma_N}}
	x^2f_{ \tilde{X_1}, \tilde{S} }(x,0)dx
	= 
	\int_{-\frac{a}{\sigma_N}}^{ \frac{1-a
		}{\sigma_N}} x^2
	f_{
		\tilde{\cal{Z}} _1,  \Sigma
	}(x,0)dx
	+O\left(N^{-1/6}{\mathbb{E}}
	\tilde{\cal{Z}}_1^2\right) +o
	\left( e^{-(\log N)^{3/2}}\right),
	\end{equation}
	and
	\begin{equation}\label{CorCov}
	\int_{-\frac{a}{\sigma_N}}^{ \frac{1-a
		}{\sigma_N}} \int_{-\frac{a}{\sigma_N}}^{ \frac{1-a
		}{\sigma_N}}xy
	f_{
		\tilde{X_j},\tilde{X_k}, \tilde{S}
	}(x,y,0)dxdy
	\end{equation}
	\[
	= 
	\int_{-\frac{a}{\sigma_N}}^{ \frac{1-a
		}{\sigma_N}} \int_{-\frac{a}{\sigma_N}}^{ \frac{1-a
		}{\sigma_N}} xy
	f_{
		\tilde{\cal{Z}}_j, \tilde{\cal{Z}}_k, \Sigma
	}(x,y,0)dxdy
	+
	O\left(N^{-1/6}{\mathbb{E}}
	\tilde{\cal{Z}}_1^2\right)
	+o
	\left( e^{-(\log N)^{3/2}}\right).
	\]
\end{cor}

\noindent
{\bf Proof.}
Let $\bar{ \cal Z} \sim{\cal N}(\bar{0}, \Lambda)$ and $\bar{X}$ be defined on the same probability space with independent distributions.
Then we 
write
\begin{equation}\label{CorV7}
f_{\tilde{X_1}, \tilde{S}}(x,0)= f_{\tilde{{\cal{Z}}}_1 , \Sigma}(x,0)
\end{equation}
\[
+ \left(f_{\tilde{S} \mid \tilde{X}_1=x}(0)
-f_{\Sigma \mid \tilde{\cal{Z}}_1=x}(0)
\right)f_{\tilde{X}_1}(x)
+\left(f_{\tilde{X}_1
}(x)-f_{\tilde{{\cal{Z}}}_1}(x)
\right)f_{\Sigma \mid \tilde{\cal{Z}}_1=x}(0).
\]
Here 
vector $\left( \tilde{{\cal{Z}}}_1,\Sigma \right) $ by its definition has a Normal distribution with zero mean and covariance matrix
\begin{equation}\label{Cm2}
\left(
\begin{array}{cc}
\frac{\Lambda_{11} }{\sigma_N^2} &  \frac{1}{N} \\ \\
\frac{1}{N}  & 1
\end{array}
\right)
,
\end{equation}
where $\Lambda_{11}\sim N^{-3}$, and $\sigma_N^2\sim N^{-2}$. Hence, the conditional distribution of 
$\Sigma$ given $\tilde{\cal{Z}}_1=x$ is again Gaussian:
\begin{equation}\label{J10*}
\Sigma \mid _{\tilde{\cal{Z}}_1=x} \ \sim \ {\cal N}\left(x\frac{\sigma_N^2}{\Lambda_{11} N} , 1- \frac{\sigma_N^2}{\Lambda_{11} N^2} \right).
\end{equation}
The variance of this conditional distribution is
\begin{equation}\label{J10}
1- O\left(\frac{1}{N}\right),
\end{equation}
implying that the density $f_{\Sigma \mid \tilde{\cal{Z}}_1=x}(0)$ as a function of $x$ is continuous and bounded uniformly in $N$ and $x$.
The latter allows us to use Theorem \ref{Lexp0} and derive
\begin{equation}\label{J11}
\int_{-\frac{a}{\sigma_N}}^{ \frac{1-a
	}{\sigma_N}}x^2
\left(f_{\tilde{X}_1
}(x)-f_{\tilde{{\cal{Z}}}_1}(x)
\right)f_{\Sigma \mid \tilde{\cal{Z}}_1=x}(0)
dx=o
\left( e^{-(\log N)^{3/2}}\right).
\end{equation}

Combining results of Theorem \ref{Lexp0} and
Lemma \ref{TJ1} (with adjustment for a conditional density)
we also get
\begin{equation}\label{J12}
\int_{-\frac{a}{\sigma_N}}^{ \frac{1-a
	}{\sigma_N}}x^2\left( f_{\tilde{S} \mid \tilde{X}_1=x}(0)
-f_{\Sigma \mid \tilde{\cal{Z}}_1=x}(0)\right)f_{\tilde{X}_1
}(x)
dx=O\left(N^{-1/6}\right)\int_{-\frac{a}{\sigma_N}}^{ \frac{1-a
	}{\sigma_N}}x^2f_{\tilde{X}_1
}(x)
dx
\end{equation}
\[=\left( {\mathbb{E}}\tilde{\cal{Z}}_1^2+ o
\left( e^{-(\log N)^{3/2}}\right) \right)  O\left(N^{-1/6}\right).\]
Now decomposition (\ref{CorV7})  together with  bounds (\ref{J12}) and (\ref{J11})  yield (\ref{CorVar}). 

The proof of  (\ref{CorCov}) follows essentially the same arguments, where conditioning on $\tilde{X}_1$ is replaced by conditioning on both variables $\tilde{X}_j$ and $\tilde{X}_k$.
\hfill$\Box$

\bigskip
\subsection{Proof of Theorem \ref{T1}}\label{MR}

Now we are ready to prove the main result. Making use of the last Corollary
\ref{CorVa} 
and Lemma \ref{TJ1} we get (recall that $a=1/N$) 
\begin{equation}\label{23Ja61*}
\mathbb{E}\left\{ \left({X}_{k,\lambda}^{\circ}-a\right)^2
\mid 
\sum_{k=1}^N X_{k, \lambda}^{\circ}=1\right\}	
=\sigma_N^2 \mathbb{E}
\left\{ \left(\frac{{X}_{k,\lambda}^{\circ}-a
}{\sigma_N}\right)^2
\ \left| \ 
\frac{\sum_{k=1}^N{X}_{k,\lambda}^{\circ}-1
}{\sigma_N}
=0
\right. \right\}
\end{equation}

\[	= \sigma_N^2 \frac{
	\int_{-\frac{a}{\sigma_N}}^{ \frac{1-a
		}{\sigma_N}} x^2f_{\tilde{X_1}, \tilde{S}}(x,0)
	d{x}
}{f_{\tilde{S}}(0)
}=  \sigma_N^2 \frac{
	\int_{-a/\sigma_N}^{ \frac{1-a
		}{\sigma_N}} x^2f_{\tilde{{\cal{Z}}}_1,  \Sigma}(x,0)
	d{x}+O\left(N^{-1/6}{\mathbb{E}}
	\tilde{\cal{Z}}_1^2\right) +o
	\left( e^{-(\log N)^{3/2}}\right)
}{f_{\Sigma}(0)
	+ O(N^{-1/6})	},
\]
where $f_{ \Sigma }(0)={1}/\sqrt{2\pi}$ as by the definition $\Sigma \sim {\cal{N}}(0,1)$. Hence, 
\begin{equation}\label{23Ja61}
\mathbb{E}\left\{ \left({X}_{k,\lambda}^{\circ}-a\right)^2
\mid 
\sum_{k=1}^N X_{k, \lambda}^{\circ}=1\right\}	
\end{equation}
\[= 
\frac{\sigma_N^2
	\int_{-a/\sigma_N}^{
		\frac{1-a
		}{\sigma_N}}
	x^2f_{ \tilde{{\cal{Z}}}_1,  \Sigma }(x,0)
	d{x}
}{f_{ \Sigma }(0)}\left(1
+  O\left(N^{-1/6}\right)\right) +o
\left( e^{-(\log N)^{3/2}}\right).
\]
Here we have (recall notation in Lemma \ref{TJ1})
\begin{equation}\label{23Ja17}
f_{ \tilde{{\cal{Z}}}_1,  \Sigma }(x,0)=\sigma_N^2
f_{{\cal{Z}}_1,  \sum_{k=1}^N {\cal{Z}}_k}(\sigma_Nx,0).
\end{equation}
By the definition $\bar{{\cal{Z}}} \sim {\cal N}(\bar{0}, \Lambda)$, 
where the covariance
matrix $\Lambda$ (see Corollary \ref{CorMay1})
with the choice (\ref{lambda}) of $\lambda$ has the following properties:
\begin{equation}\label{23J64}
{\bf Var}({\cal{Z}}_1)=\Lambda_{11} \sim N^{-3}, \ \ \  
{\bf Cov} \left( {\cal{Z}}_1,\sum_{k=1}^N {\cal{Z}}_k \right) =\frac{\sigma_N^2}{N} \sim N^{-3}.
\end{equation}
Then making use of the closed form for the density
in (\ref{23Ja17}), which is
\begin{equation}\label{23J63}
f_{{\cal{Z}}_1,  \sum_{k=1}^N {\cal{Z}}_k}(\sigma_Nx,0)
=\frac{1}{2\pi \sqrt{\Lambda_{11} \sigma_N^2} \sqrt{1- \frac{1}{\Lambda_{11} \sigma_N^2 }\left(\frac{\sigma_N^2}{ N} \right)^2 }}\exp
\left\{
-\frac{1}{2\left(1- 
	\frac{\sigma_N^2}{\Lambda_{11}   N^2 }\right)}
\frac{(\sigma_N x)^2}{\Lambda_{11} }
\right\},
\end{equation}
and the asymptotic (\ref{23J64}), it is straightforward to compute that 
\begin{equation}\label{23J65}
\int_{|x|>a/\sigma_N}
x^2
\sigma_N^2
f_{{\cal{Z}}_1,  \sum_{k=1}^N {\cal{Z}}_k}(\sigma_Nx,0)
d{x}=o(e^{-\sqrt{N}}).
\end{equation}
This allows us to derive from (\ref{23Ja61})
\begin{equation}\label{23Ja65}
\mathbb{E}\left\{ \left({X}_{k,\lambda}^{\circ}-a\right)^2
\mid 
\sum_{k=1}^N X_{k, \lambda}^{\circ}=1\right\}	
\end{equation}
\[= 
\frac{
	\int_{-\infty}^{+\infty}x^2
	f_{{\cal{Z}}_1,  \sum_{k=1}^N {\cal{Z}}_k}(x,0)
	d{x} 
}{f_{ \sum_{k=1}^N {\cal{Z}}_k }(0) }(1+O\left(N^{-1/6}\right))
+o(e^{-\sqrt{N}})+o
\left( e^{-(\log N)^{3/2}}\right) \]
\[={\mathbb E}
\left\{ 
{\cal{Z}}_k^2
\mid 
\sum_{k=1}^N 
{\cal{Z}}_k=0\right\}(1+O\left(N^{-1/6}\right))+o
\left( e^{-(\log N)^{3/2}}\right)
.\]

In precisely same manner we derive with a help of (\ref{CorCov}) for any $j>1$ the following (similar to (\ref{23Ja61})) relations:
\begin{equation}\label{23Ja3}
{\bf Cov}\left\{ \left( X_{1,\lambda}^{\circ}, X_{j,\lambda}^{\circ}\right) \mid \sum_{k=1}^N X_{k, \lambda}^{\circ}=1\right\}
\end{equation}
\[
=\sigma_N^2 \mathbb{E}\left\{ \left(\frac{
	{X}_{1,\lambda}^{\circ}-a
}{\sigma_N}\right) \left( \frac{
	{X}_{j,\lambda}^{\circ}-a
}{\sigma_N}\right)
\ \left| \ 
\frac{\sum_{k=1}^N{X}_{k,\lambda}^{\circ}-1
}{\sigma_N}
=0
\right. \right\}
\]
\[	= \sigma_N^2 \frac{
	\int_{-\frac{a}{\sigma_N}}^{ \frac{1-a
		}{\sigma_N}} \int_{-\frac{a}{\sigma_N}}^{ \frac{1-a
		}{\sigma_N}} 
	xyf_{\tilde{X_1}, \tilde{X_j},\tilde{S}}(x,y,0)
	d{x}dy
}{f_{\tilde{S}}(0)
}\]
\[=  \sigma_N^2 
\frac{
	\int_{-\frac{a}{\sigma_N}}^{ \frac{1-a}{\sigma_N}}
	\int_{-\frac{a}{\sigma_N}}^{ \frac{1-a}{\sigma_N}}
	xy
	f_{ \tilde{{\cal{Z}}}_1,  \tilde{{\cal{Z}}}_j,\Sigma}(x,0)
	d{x}dy+
	O\left(N^{-1/6}{\mathbb{E}}
	\tilde{\cal{Z}}_1^2\right)
	+o
	\left( e^{-(\log N)^{3/2}}\right)
}{
	f_{\Sigma}(0)
}
.
\]

It is again straightforward to derive from here with a help of Proposition \ref{Pcov} that 
\begin{equation}\label{23Ja70}
{\bf Cov}\left\{ \left( X_{1,\lambda}^{\circ}, X_{j,\lambda}^{\circ}\right) \mid \sum_{k=1}^N X_{k, \lambda}^{\circ}=1\right\}
\end{equation}
\[
=\sigma_N^5 
\frac{
	\int_{-\infty}^{ +\infty}
	\int_{-\infty}^{+\infty}
	xy
	f_{ {{\cal{Z}}}_1,  {{\cal{Z}}}_j,\sum_{k=1}^N {\cal{Z}}_k}(\sigma_Nx,\sigma_Ny, 0)
	d{x}dy
}{
	f_{\Sigma}(0)
}
+ O\left(\Lambda _{11}N^{-1/6}\right)
\]
\[
=
\frac{
	\int_{-\infty}^{ +\infty}
	\int_{-\infty}^{+\infty}
	xy
	f_{ {{\cal{Z}}}_1,  {{\cal{Z}}}_j,\sum_{k=1}^N {\cal{Z}}_k}(x,y, 0)
	d{x}dy
}{
	f_{ \sum_{k=1}^N {\cal{Z}}_k}(0)
}
+  O\left(\Lambda _{11}N^{-1/6}\right)
\]
\[={\bf Cov}
\left\{ 
\left( {\cal{Z}}_1,  {\cal{Z}}_j\right)
\mid 
\sum_{k=1}^N 
{\cal{Z}}_k=0\right\}
+  O\left(N^{-3-1/6}\right)
.\]
By the properties of normal distribution the conditional distribution of
\[
\left( {\cal{Z}}_1,  {\cal{Z}}_j\right)
\mid 
\sum_{k=1}^N 
{\cal{Z}}_k=0 
\]
is again normal with zero mean vector and the covariance matrix
\begin{equation}\label{23J76}
\left(
\begin{array}{cc}
\Lambda_{11} & \Lambda_{1j}  \\ \\
\Lambda_{1j} & \Lambda_{jj}  
\end{array}
\right)
-\frac{\sigma_N^2}{N^2}\left(
\begin{array}{cc}
1 & 1 \\ \\
1 & 1
\end{array}
\right).
\end{equation}
Observe that by the definition
\[\sigma_N^2=N\left( \Lambda_{11} + 2\sum_{k=2}^{\left[ \frac{N}{2}\right]} \Lambda_{1j}\right).\]
Hence, applying the formula from Proposition \ref{Pcov} with $\lambda$ chosen as in (\ref{lambda}) we derive from  (\ref{23J76})
\[
{\bf{Cov}}\left\{
\left( {\cal{Z}}_1,  {\cal{Z}}_j\right)
\mid 
\sum_{k=1}^N 
{\cal{Z}}_k=0 
\right\}=\Lambda_{1j}-\frac{\sigma_N^2}{N^2}=\Lambda_{11} \left(
(-\delta)^{j-1}-\frac{1}{N}\frac{1-\delta}{1+\delta} +o\left(\delta ^{N/2} \right)
\right)
\]
\[
= \frac{1}{N^3}
\frac{1}{2\beta + \gamma (1 - \delta)/2}
\left(
(-\delta)^{j-1}-\frac{1}{N}\frac{1-\delta}{1+\delta} +o\left(e ^{-\alpha}N \right)
\right)
\]
for some positive constant $\alpha$.
This together with  (\ref{23Ja70}) and (\ref{23Ja65}) confirms the statements of Theorem \ref{T1}. \hfill$\Box$
\bigskip

\noindent
{\bf Proof of Theorem \ref{T2o}.}
To proof closely follows the lines of above proof replacing results of 
Theorem \ref{Lexp0} and Corollary \ref{ExCov} by the results of 
Theorem \ref{Lexp1} and Corollary \ref{ExCov1}, correspondingly. Then we conclude that the conditional distribution $(X_{1,\lambda}, X_{j,\lambda}) \mid \sum_{k=1}^N X_{k, \lambda}=1$ is related to the distribution of ${\bar X}$ described by ${\cal Z}'$ in
Theorem \ref{Lexp1}
in precisely same manner as for the circular case.

This formal argument enables us to use approximation (\ref{Ea2}) but now indirectly for ${\cal H}_1^{-1}$, whose entries are proved to be close to the covariances, and the result of Theorem \ref{T2o} follows. \hfill$\Box$

\section*{Appendix} 
\subsection*{Inverse of a circular matrix}
The inverses of circular matrices and closely related Toeplitz matrices often arise in the study of models of mathematical physics (\cite{DIK}). Some direct computational results are also available (e.g., \cite{C}, \cite{K}), however typically they are very limited to particular cases. As here we are interested only in certain asymptotic rather than in deriving all entries for the inverse matrices, we can avoid heavy computations looking  instead for the asymptotic directly.

For any $n\geq 4$ let 
${\cal H}(n)$ be $n \times n$  matrix defined (as in (\ref{Hess})) by 
\begin{equation}\label{Hessn}
{\cal H}(n):= 
\frac{1}{a^3}\left(
\begin{array}{cccccc}
2\beta+\frac{\gamma}{2} & \frac{\gamma}{4} & 0 & \ldots&  0 & \frac{\gamma}{4} \\
\frac{\gamma}{4} & 2\beta+\frac{\gamma}{2}& \frac{\gamma}{4} & 0&  \ldots&  0 \\
\ldots & & & & & \\
\frac{\gamma}{4} & 0 &  \ldots &   0 & \frac{\gamma}{4} & 2\beta+\frac{\gamma}{2}
\end{array}
\right).
\end{equation}
Note that here $a=a(N)$, and we treat $n$ as an independent variable.
Consider
\begin{equation}\label{May1^*}
\Lambda (n):=  {\cal H}^{-1}(n)=\left(\Lambda_{ij}(n)\right)_{1\leq i,j\leq N}.
\end{equation}
First we note, that
by the symmetry we have in our case  for all $n\geq 4$ 
\begin{equation}\label{syminv}
({\cal H}^{-1})_{ii}=({\cal H}^{-1})_{11}=\Lambda_{11},
\end{equation}
as well as
\begin{equation}\label{syminv1}
({\cal H}^{-1})_{ij}=({\cal H}^{-1})_{1 (1+|i-j|)}= \Lambda_{1(1+|j-i|)},
\end{equation}
for all $1\leq i,j\leq N.$
Let us compute now
\begin{equation}\label{inv1}
\Lambda_{ii}=({\cal H}^{-1})_{11}= \frac{\det {\cal D}_{n-1}}{\det {\cal H}(n)} ,
\end{equation}
where ${\cal D}_{n-1}$ is the 3-diagonal  symmetric $(n-1)\times (n-1)$ matrix corresponding to ${\cal H}(n)$:
\begin{equation}\label{Diag}
{\cal D}_{n-1}:= 
\frac{1}{a^3}\left(
\begin{array}{cccccc}
2\beta+\frac{\gamma}{2} & \frac{\gamma}{4} & 0 & \ldots&  0 & 0 \\
\frac{\gamma}{4} & 2\beta+\frac{\gamma}{2}& \frac{\gamma}{4} & 0&  \ldots&  0 \\
\ldots & & & & & \\
0 & 0 &  \ldots &   0 & \frac{\gamma}{4} & 2\beta+\frac{\gamma}{2}
\end{array}
\right).
\end{equation}

To simplify notation define for any $n>4$ and $b,c \in \mathbb{R}$ an $n\times n$ circular matrix
\begin{equation}\label{May1}
M_n:= 
\left(
\begin{array}{cccccc}
c  & b  & 0 & \ldots&  0 & b \\
b & c & b & 0&  \ldots&  0 \\
\ldots & & & & & \\
0  & 0  & \ldots &  b &  c & b  \\
b & 0 &  \ldots &   0 & b & c
\end{array}
\right).
\end{equation}
Then $M_n={\cal H}(n)$ if
\begin{equation}\label{May2}
b= \frac{1}{a^3} \frac{\gamma}{4}, \ \ \ c=\frac{1}{a^3}\left( 2\beta+\frac{\gamma}{2} \right) = \frac{2\beta}{a^3} + 2 b.
\end{equation}
The eigenvalues of $M_n$ by (\ref{eigenM}) are
\begin{equation}\label{eigenMv}
\nu^{\circ}_j=c +2b\cos \left( 2\pi \frac{ j}{n}\right), \ \ \ j=1, \ldots, n.
\end{equation}
Correspondingly, we rewrite
${\cal D}_{n}$ introduced in (\ref{Diag})
as
\begin{equation}\label{MayDiag}
{\cal D}_n = 
\left(
\begin{array}{cccccc}
c & b & 0 & \ldots&  0 & 0 \\
b & c & b & 0&  \ldots&  0 \\
\ldots & & & & & \\
0 & 0 &  \ldots &   0 & b & c
\end{array}
\right).
\end{equation}
The inverse of $M_n$ we shall also denote $\Lambda$:
\[
M_n^{-1}=\Lambda = (\Lambda_{ij})_{1\leq i,j\leq n}.
\]
Exploring the symmetry of the introduced matrices it is straightforward to derive the following relations (one may also consult \cite{D} or \cite{C}): 
\begin{equation}\label{MayInv}
\begin{array}{l}
\Lambda_{11} = \frac{1}{{\mbox {\bf det}} M_n} D_{n-1} ,\\ \\ 
\Lambda_{12}=\Lambda_{1n} = \frac{1}{{\mbox {\bf det}} M_n}\left(
-bD_{n-2}+(-b)^{n-1}\right),\\ \\ 
\Lambda_{1k} = \frac{(-1)^{k+1}}{{\mbox {\bf det}} M_n}\left(  b^{k-1}D_{n-k}+(-1)^{n-2}b^{n-k+1}D_{k-2}
\right)
, \ \ k =2, \ldots, n-1.
\end{array}
\end{equation}
The eigenvalues of ${\cal D}_n$  are also known (see, e.g., \cite{K}), these are
\begin{equation}\label{3eigenM}
\nu_j=c + 2b\cos \left( \pi \frac{ j}{n+1}\right), \ \ \ j=1, \ldots, n.
\end{equation}
Hence, we have  the determinant
\begin{equation}\label{det M}
D_n:={\mbox {\bf det} \ } {\cal D}_n= \prod_{j=1}^n\left( c + 2b\cos \left( \pi \frac{ j}{n+1}\right)\right).
\end{equation}
In our case $c>2b $
by (\ref{May2}), and therefore we derive from the last formula
\[
\frac{1}{n} \log D_n= \log c 
+\frac{1}{n} \sum_{j=1}^n \log 
\left( 1 + \frac{2b}{c}\cos \left( \pi \frac{ j}{n+1}\right)\right) ,
\]
which converges to 
\[
\log c + \int_0^1 \log \left( 1 + \frac{2b}{c}\cos \left( \pi x\right)\right) dx= \log c + \frac{1}{\pi }\int_0^{\pi} \log \left( 1 + \frac{2b}{c}\cos x\right) dx,
\]
as $n \rightarrow \infty.$ 
Computing the last integral  we get from here
\begin{equation}\label{May10}
\lim_{n \rightarrow \infty}\frac{1}{n} \log D_n = \log c + \log \frac {1 + \sqrt{ 1-\left(\frac{2b}{c}\right)^2}}{2} = \log \frac {c + \sqrt{ c^2-(2b)^2}}{2}.
\end{equation}

Observe, that the last asymptotic can be alternatively derived even without use of eigenvalues. This follows simply by the recurrent relation for the determinant 
\begin{equation}\label{May11}
D_n:={\mbox {\bf det} \ } {\cal D}_n=cD_{n-1}-b^2D_{n-2},
\end{equation}
whose characteristic equation
\begin{equation}\label{May12}
x^2-
cx+b^2=0
\end{equation}
has the roots 
\begin{equation}\label{May13}
x_{1} = \frac{c+ \sqrt{ c^2-4b^2}}{2}, \ \ \ \ 
x_{2} = \frac{c- \sqrt{ c^2-4b^2}}{2}
.
\end{equation}
Hence, 
\begin{equation}\label{May14}
D_n= A
x_{1}^n + Bx_{2}^n
\end{equation}
where 
taking into account that 
\[D_1=c, \ \ D_2= c^2-b^2,\]
the constants $A,B$ satisfy
\[
\left\{
\begin{array}{ll}
A
x_{1} + Bx_{2}& =c,\\ \\
A
x_{1}^2 + Bx_{2}^2& =c^2-b^2.
\end{array}
\right.
\]
As in our case $c>2b>0$ this gives us a non-zero value $A$, hence the largest root $x_1$ 
gives us the same asymptotic (\ref{May10}).

To find the remaining determinant ${\mbox {\bf det}} M_n$ in 
(\ref
{MayInv}) one could also use the eigenvalues (\ref{eigenM}), but instead in the following we shall explore the relation
\begin{equation}\label{May15}
{\mbox {\bf det}} M_n = cD_{n-1}-2b^2D_{n-2}-2(-1)^{n}b^n=D_{n}-b^2D_{n-2}-2(-1)^{n}b^n,
\end{equation}
together with observation that by (\ref{May13}) and (\ref{May2})
\begin{equation}\label{May16}
x_{1} = \frac{c+ \sqrt{ c^2-4b^2}}{2}
>\frac{c}{2}>b.
\end{equation}

Now we ready to derive the asymptotic for the terms in (\ref{MayInv}).

Combining (\ref{May15}) and (\ref
{MayInv}) we get
\begin{equation}\label{May22}
\Lambda _{11}=\frac{D_{n-1}
}{cD_{n-1}-2b^2D_{n-2}-2(-1)^{n}b^n}
\end{equation}
\[=\frac{x_1
}{cx_{1}-2b^2} +O(e^{-\alpha n})
= \frac{1}{c-2b\frac{
		b}{x_1}}+O(e^{-\alpha n}),
\]
where $\alpha>0$ is some positive constant.

In a similar way by (\ref{May15}) and (\ref
{MayInv})  we have
\begin{equation}\label{May23}
\Lambda _{12}=\Lambda _{1n}=\frac{-bD_{n-2}+(-b)^{n-1}
}{cD_{n-1}-2b^2D_{n-2}-2(-1)^{n}b^n}
\end{equation}
\[=-\frac{x_1
}{cx_{1}-2b^2} \frac{b}{x_1}+O(e^{-\alpha n}), \]
and also for all $2<k\leq \left[
\frac{n+1}{2} \right]$
\begin{equation}\label{May24}
\Lambda_{1k} = \Lambda_{1 (n-k+2)} =\frac{
	(-b)^{k-1}D_{n-k}+(-1)^{n+k+1}b^{n-k+1}D_{k-2}
}{cD_{n-1}-2b^2D_{n-2}-2(-1)^{n}b^n}
\end{equation}

\[=
\frac{
	(-b)^{k-1}x_1^{2-k}+(-b)^{n-k+1}x_1^{k-n}
}{cx_1-2b^2}+O(e^{-\alpha n})\]

\[=(-1)^{k-1}
\left(\frac{
	b}{x_1}\right)^{k-1}
\frac{x_1}{cx_1-2b^2}+O(e^{-\alpha n/2}).\]

Rewriting the last formulae in terms of $\delta, \beta , \gamma $ with a help of (\ref{May2}) and (\ref{May13}) we arrive at the statement of Proposition \ref{Pcov}. \hfill$\Box$

\end{document}